\documentclass[lettersize,journal]{IEEEtran}
\usepackage{amsmath,amsfonts}
\usepackage{algorithmic}
\usepackage{algorithm}
\usepackage{array}
\usepackage{textcomp}
\usepackage{stfloats}
\usepackage{url}
\usepackage{verbatim}    
\usepackage{graphicx}
\usepackage{amssymb}
\usepackage{hyperref}                   
\usepackage{subfigure}     
\usepackage{fancyhdr}      

\makeatletter
\renewcommand{\maketag@@@}[1]{\hbox{\m@th\normalsize\normalfont#1}}%
\makeatother

\newtheorem{theorem}{Theorem}
\newtheorem{lemma}{Lemma}
\newtheorem{remark}{Remark}

\begin{document}

\title{Delay Performance Analysis of Delay-Deterministic Wireless Networks with Infinite and Finite Blocklength Transmission}

\author{Hanxue~Ding,  
	Shaoyi~Xu, \textit{Member}, \textit{IEEE},
	Ziheng~Xu,
	Rongtao~Xu,
	Zonghui~Li,
	and Junhui~Zhao, \textit{Senior Member}, \textit{IEEE}

\thanks{Hanxue Ding, Shaoyi Xu, Ziheng Xu, Rongtao Xu and Junhui Zhao are with the School of Electronic and Information Engineering, Beijing Jiaotong University, Beijing 100044, China (e-mail: hanxue\_ding@bjtu.edu.cn; shyxu@bjtu.edu.cn; 22115012@bjtu.edu.cn; rtxu@bjtu.edu.cn; junhuizhao@bjtu.edu.cn).
	
Zonghui Li is with the School of Computer and Information Technology, Beijing Jiaotong University, Beijing 100044, China (e-mail: lizonghui@bjtu.edu.cn).}
}

\pagestyle{fancy}
\fancyhf{}
\fancyhead[C]{\fontsize{9}{9}\selectfont This work has been submitted to the IEEE for possible publication. Copyright may be transferred without notice, after which this version may no longer be accessible.}
\renewcommand{\headrule}{}



\maketitle

\begin{abstract}
In order to achieve stable and reliable industrial manufacturing, wireless networks must meet the stringent communication requirements of industrial automation, particularly the need for deterministic low latency communication. The limited wireless resources and time-varying fading channel contribute to the random fluctuations of transmission delay, making it challenging to realize delay-deterministic wireless networks. An open challenge in this context is to model delay determinism, also known as jitter, and analyze delay performance. In this paper, we model jitter as the variance of delay and conduct a comprehensive analysis of delay performance. Specifically, we consider two transmission regimes: infinite blocklength (IBL) and finite blocklength (FBL). In the IBL regime, the distribution of the transmission delay is analyzed, and the closed-form expressions for the average delay, jitter, and delay violation probability are derived. In the FBL regime, an upper bound on the transmission delay is first approximated at a high signal-to-noise ratio. Based on this upper bound, the delay distribution, delay violation probability, average delay, and jitter are derived. Finally, simulation results are provided to validate the accuracy of the analysis and derivations. Additionally, the impact of system parameters on jitter is analyzed to gain further insights.
\end{abstract}


\section{Introduction} \label{sec1}
\IEEEPARstart{W}{ith} the deep integration of wireless communication technology and traditional industry, the Industrial Internet of Things (IIoT), as a new application paradigm, has attracted extensive attention from academia and industry. IIoT is considered an important cornerstone of Industry 4.0 for its advantages of low cost, easy deployment, easy maintenance, and high flexibility \cite{Vitturi}. Typical IIoT applications need to realize delay-deterministic wireless networks, which need to meet stringent communication requirements, including high reliability, strong real-time, low latency, and determinism \cite{Luvisotto}.

Transmission delay is an inherent characteristic and a crucial performance metric of communication systems. Traditional industrial networks mostly adopt wired infrastructure, which provides a stable and reliable transmission process with a fixed delay. Due to the open and random nature of wireless channels, the transmission process is unpredictable and unreliable, resulting in random fluctuations in transmission delay. A significant number of metal surfaces in the industrial environment cause multiple signal reflections, making the industrial wireless channel usually a severe multipath fading channel \cite{Liu}. Thus, the fluctuations in delay are more severe. To this end, the main challenge of replacing wired connections with wireless transmission in industrial networks is the random transmission delay. For this reason, industrial wireless networks require deterministic low delay.

To enable low-latency communication, finite blocklength coding is considered at the physical layer to reduce the transmission delay by reducing blocklength. Different from Shannon capacity, which assumes error-free transmission, the decoding error probability is non-negligible due to finite blocklength coding. Taking into account the blocklength and decoding error probability, the maximum achievable coding rate in the finite blocklength (FBL) regime was recharacterized in \cite{Polyanskiy}. Inspired by \cite{Polyanskiy}, low-latency communication based on finite blocklength coding has been extensively researched. On the one hand, resource allocation and scheduling algorithms based on finite blocklength coding are commonly studied to ensure low-latency and high-reliability communication \cite{Ren}, \cite{Lin}, \cite{Xu}. The authors in \cite{Ren} jointly optimized the blocklength and transmit power of a short packet communication system to maximize reliability while satisfying the maximum delay constraint. The authors in \cite{Lin} achieved the minimum transmission delay while guaranteeing the reliability constraint through the joint design of channel training and data transmission. A packet scheduling strategy based on a hybrid non-orthogonal multiple access (NOMA) scheme has been proposed to serve users with heterogeneous delay requirements \cite{Xu}. On the other hand, many researchers have attempted to reduce delay by introducing other techniques into finite blocklength transmission \cite{Sun}, \cite{Yilmaz}, \cite{Ghanem}, \cite{Hashemi}. For example, the authors in \cite{Sun} demonstrated the superiority of NOMA over orthogonal multiple access in reducing transmission delay while achieving the same effective throughput. Reference \cite{Yilmaz} revealed that using multiple antennas for transmit or receive diversity can reduce transmission delay. The study in \cite{Ghanem} examined the significance of beamforming in reducing the latency and enhancing the reliability of ultra-reliable and low-latency communication (URLLC) systems. The authors in \cite{Hashemi} illustrated the important role of RIS in reducing transmission delay. However, the above works did not consider deterministic delay.

Deterministic delay imposes a new requirement on the stability of the delay, i.e., low jitter. There have been several studies on deterministic delay. The works \cite{Chen} and \cite{H. Xu} defined delay determinism as the probability that delay falls within a specified time window and maximized delay determinism by jointly optimizing power and bandwidth. In \cite{L. Chen}, the authors characterized jitter as the difference between the maximum and average delay and jointly allocated computation and communication resources to minimize end-to-end jitter. Strict constraints were imposed on the upper and lower bounds of the delay to reduce jitter \cite{Li}. Penalties were applied to packets that arrived either too early or too late to ensure delay determinism \cite{Sheikh}. The above studies \cite{Chen}, \cite{H. Xu}, \cite{L. Chen}, \cite{Li}, \cite{Sheikh} have ensured clear upper and lower bounds on latency in various ways. In addition to the above studies, some other works define jitter as the variance of the delay \cite{Wei}, \cite{Liu}. Specifically, the authors in \cite{Wei} pointed out that delay-deterministic wireless networks need to support low-delay and low-jitter data transmission, proposed spectrum sharing between industrial and cellular networks to improve the delay performance of industrial networks, and derived the jitter with and without spectrum sharing. The simulation results illustrated that spectrum sharing can effectively reduce delay and jitter. To cope with the issue of industrial control instability caused by random transmission delay and jitter, a two-layer online control algorithm was designed to optimize the delay distribution under the fading channel and minimize delay variance \cite{Liu}. Although considering delay determinism, the study in \cite{Liu} made an over-optimistic assumption about infinite transmission power, which is unrealistic for real-world systems. Note that the above works are developed based on the classical Shannon capacity formula and are not applicable to finite blocklength transmission scenarios.

In essence, deterministic low latency is the combination of low latency and low jitter. Although the above-mentioned works have investigated low-latency and low-jitter transmission separately, they have not considered other delay properties. Existing studies on delay performance analysis mainly include three aspects. \textit{1) Mean and variance of delay.} Variance is a statistical measure that represents the degree of dispersion of delay and is therefore often used to characterize jitter, as mentioned in \cite{Liu} and \cite{Wei}. \textit{2) Statistical distribution.} The statistical distribution of random variables contains rich statistical information, such as mean and variance. In addition, delay distribution has profound implications for resource allocation under delay constraints. A unified approximation framework for the delay distribution of federated learning (FL) was presented in \cite{L. Li}. To be specific, the distribution of the overall delay was derived based on the saddle-point approximation, and the tail distribution was acquired based on the extreme value theory and large deviation theory. Because of the large number of parameters in the FL, the delay distribution in \cite{L. Li} was derived based on Shannon capacity. Reference \cite{Danufane} considered discrete slots and ensured that the data could be successfully received by retransmission, so the transmission delay was a discrete random variable related to the slot number. However, this paper focuses on a one-shot transmission, and the transmission delay is a continuous variable. \textit{3) Probabilistic delay bounds.} An upper bound on delay violation probability is described based on queuing theory. For example, the works \cite{Schiessl}, \cite{S. Schiessl}, \cite{P. Cui} analyzed the transmission delay at the physical layer and the queuing delay at the upper layer and derived the delay violation probability based on stochastic network calculus. While the end-to-end delay is usually studied in the above works, this paper only focuses on the transmission delay.\footnote[1]{On the one hand, the primary difficulty in replacing wired connections with wireless transmission is that wireless transmission delay is random, so we focus on transmission delay. On the other hand, queuing delay is usually analyzed using queuing theory and stochastic network calculus, which significantly differ from the methods used to analyze transmission delay. Therefore, the queuing delay is not currently considered in this paper.} In spite of the fact that many published works have separately investigated low delay, delay determinism, and different delay performances, a comprehensive analysis of transmission delay performance in the FBL regime is still an open problem. 

Against this background, aiming to guarantee the deterministic low delay requirement, we propose a unified framework to analyze the transmission delay performance, including the statistical distribution, delay violation probability, average delay, and jitter. To the best of our knowledge, this is the first work that provides a comprehensive analysis of delay performance in the FBL regime. The main contributions are summarized as follows.

\begin{itemize}
	\item In the infinite blocklength (IBL) regime, we first analyze the cumulative distribution function (CDF) and the probability density function (PDF) of the transmission rate and transmission delay for the single-antenna case. Then, we show that the single-antenna scenario is a special case of the multi-antenna scenario. Consequently, the following analysis is based on multi-antenna transmission. We also derive the closed-form expressions for the average delay, jitter, and delay violation probability.
	
	\item In the FBL regime, the transmission delay and theoretical distribution of transmission delay are first provided. Since the expression for delay is too complicated to support subsequent analysis, an approximation is obtained at a high signal-to-noise ratio (SNR). From this approximation, the approximate delay  distribution and delay violation probability are derived, and the closed-form expressions for average delay and jitter are obtained. 
	
	\item Extensive simulations demonstrate that the theoretical and approximation results are in good agreement with the simulation results. We compare the delay performance of single-antenna and multi-antenna scenarios and finite and infinite blocklength transmissions. Based on these comparisons, we draw useful insights into the design and optimization of delay-deterministic wireless networks. Moreover, related parameters such as average SNR and the number of transmit antennas should be carefully considered to improve delay performance.	
\end{itemize}

The rest of this paper is organized as follows. Section \ref{sec2} presents the system model, which contains both single-antenna and multi-antenna scenarios. Section \ref{sec3} analyzes the delay performance of long packet communication based on Shannon capacity. Section \ref{sec4} analyzes the delay performance of short packet communication. Simulation results and analysis are presented in Section \ref{sec5}. Finally, Section \ref{sec6} concludes the paper.

\textit{Notations:} 
In this paper, $\mathbf{I}_N$ denotes a size-$N$ identity matrix. The operators $\mathbb{E}[\cdot] $ and $\text{Var}[\cdot]$ denote the expectation and variance, respectively. $X \sim \mathcal{CN}(\mu_X, \sigma^2_X)$ represents the complex Gaussian variable with mean $\mu_X$ and variance $\sigma^2_X$. $X \sim \text{Gamma}(a, c)$ represents the gamma random variable with shape and rate parameters $a$ and $c$, respectively. $X \sim \text{Exp}(\lambda)$ denotes the exponential random variable with rate parameter $\lambda$. $F_X(\cdot)$ and $f_X(\cdot)$ denote the CDF and PDF of the random variable $X$, respectively. 

\section{System Model} \label{sec2}

We consider a wireless control network in an IIoT scenario where a central controller serves a single antenna actuator. The controller transmits control messages to the actuator. After receiving the control messages, the actuator takes action based on the received messages. The control messages contain $L$ bits of information and are encoded into $n$ symbols. Due to the limited wireless resources and the time-varying fading channel, the transmission delay randomly fluctuate in a large range, resulting in a large jitter. Compared to maximum transmission delay, jitter is more challenging for system stability \cite{Luvisotto}. Therefore, jitter is our primary performance metric of interest. Similar to references \cite{Liu} and \cite{Wei}, we define jitter as delay variance.

Assume that the large-scale fading coefficient between the controller and the actuator is given by $\beta = \chi_0 d^{-\alpha}$, where $\chi_0$ represents the power gain at the reference distance, $ d $ denotes the distance between the controller and the actuator, and $ \alpha $ is the path loss factor. The small-scale channel $h$ is assumed to experience Rayleigh fading.\footnote[2]{The measurement results demonstrate that indoor industrial environments have different channel distributions at different locations, such as Rayleigh, Rician, and Nakagami-$m$ distributions \cite{Cheffena}. For the sake of simplicity and clarity, we assume that the channel follows Rayleigh distribution, and the same assumption is also used in \cite{Ren}, \cite{Yilmaz}, \cite{Wei}. However, the analytical framework of this paper can be extended to other channel models.} A quasi-static fading channel is considered, which implies that the wireless channel remains constant in the transmission duration of a packet. 

If the controller has a single antenna, i.e., $ h \sim \mathcal{CN} (0,1) $, the received SNR can be represented as
\begin{equation}
	\label{eq1}
	\gamma = \frac{P_t \beta |h|^2}{\sigma^2} \triangleq \rho |h|^2 , \tag{1}  
\end{equation}
where $ \rho \triangleq \frac{P_t \beta}{\sigma^2} $ is the average SNR, $ P_t $ is the transmit power, and $ \sigma^2 $ is the noise power. The received SNR follows an exponential distribution, that is, $ \gamma \sim \text{Exp} \left(\frac{1}{\rho}\right) $. The CDF and PDF of $\gamma$ are given by
\begin{equation*} 
	F_{\gamma}(x)= 1 - e^{-\frac{1}{\rho} x}, 
	\tag{2}   \label{eq2} 
\end{equation*}
\begin{equation*} 
	f_{\gamma}(x) = \frac{1}{\rho} e^{-\frac{1}{\rho}x} . \tag{3}  \label{eq3}
\end{equation*}

If the controller has $N$ antennas, i.e., $ \mathbf{h} \sim \mathcal{CN} (\mathbf{0},\mathbf{I}_N) $, the received SNR is
\begin{equation*} \gamma = \frac{P_t \beta \left \|  \mathbf{h} \right \|^2}{N\sigma^2} = \frac{\rho}{N} \left \| \mathbf{h} \right \|^2 , \tag{4}  \label{eq4}
\end{equation*}
where the received SNR follows a gamma distribution, that is, $ \gamma \sim \text{Gamma}(N, \frac{\rho}{N}) $. The CDF and PDF of $\gamma$ are formulated as
\begin{equation*} 
	F_{\gamma}(x)= \frac{\gamma(N,\frac{N}{\rho}x)}{\Gamma(N)} , \tag{5}  \label{eq5}
\end{equation*}
\begin{equation*} 
	f_{\gamma}(x) = \frac{N^N x^{N-1}  e^{-\frac{N}{\rho}x }}{\rho^{N}\Gamma(N)},  \tag{6}  \label{eq6}
\end{equation*}
where $ \gamma(\cdot, \cdot) $ is the lower incomplete gamma function \cite[eq. (8.350.1)]{Gradshteyn} and $ \Gamma(N) = (N-1)! $.

\section{Infinite Blocklength Analysis} \label{sec3}
The IBL assumption is idealistic for this scenario. However, for the sake of completeness of the analysis and to compare the corresponding results of IBL and FBL, we first derive delay performance in the IBL regime.

In the IBL regime, the maximum coding rate for error-free transmission is given by Shannon capacity, which is
\begin{equation*} 
	R = \log_2 (1+\gamma). \tag{7}  \label{eq7}
\end{equation*}
Assuming that the transmission bandwidth is $B$ Hz, then the transmission delay is denoted by
\begin{equation*} 
	\tau = \frac{L}{BR} = \frac{L}{B\log_2 (1+\gamma) }. \tag{8}  \label{eq8}
\end{equation*}
It can be seen that the transmission delay is related to the number of transmitted bits, bandwidth, and instantaneous SNR. Jitter is defined as the variance of $\tau$, which is
\begin{equation*} 
	\text{Var} [\tau ] = \mathbb{E} [\tau^2] - \mathbb{E} [\tau]^2. \tag{9}  \label{eq9}
\end{equation*}
According to the definition of expectation, the moment of delay can be expressed as
\begin{equation*} 
	\mathbb{E} [\tau^i ] = \int_0^{\infty} t^i f_{\tau} (t) dt. \tag{10}  \label{eq10}
\end{equation*}
As seen from \eqref{eq10}, the closed-form expression for jitter requires knowledge of the PDF of the transmission delay. In what follows, we will analyze the distribution of delay in single-antenna and multi-antenna scenarios separately and then calculate the jitter.

\subsection{Delay Distribution in the Single-Antenna Scenario} \label{sec3.1}
In the single-antenna scenario, we derive the CDF and PDF of the transmission rate and transmission delay based on \eqref{eq2} and \eqref{eq7}, which are given by Lemma \ref{lemma1} and Lemma \ref{lemma2}, respectively.

\begin{lemma}  \label{lemma1}
	The CDF and PDF of the transmission rate in the single-antenna scenario are formulated as
	\begin{equation*} 
		F_{R}(y) = 1 - e^{-\frac{1}{\rho}\left(2^y-1\right)}, \tag{11} \label{eq11}
	\end{equation*}
	\begin{equation*} 
		f_{R}(y) = \frac{\ln 2}{\rho} 2^y e^{-\frac{1}{\rho} \left( 2^y - 1 \right)}. \tag{12} \label{eq12}
	\end{equation*}
\end{lemma}

\textit{Proof:}
According to the expression of the transmission rate and the definition of CDF, the CDF of the transmission rate can be expressed as
\begin{align*} 
	F_{R}(y) &= \Pr \left\{ \log_2(1+ \gamma) \le y \right\}   \\
	& = \Pr \left\{ \gamma \le 2^y - 1 \right\}   \\
	& = F_{\gamma} \left( 2^y-1 \right). \tag{13} \label{eq13}
\end{align*}
Substituting \eqref{eq2} into the above equation yields the CDF of the transmission rate in \eqref{eq11}. The derivation of the CDF yields the PDF of the transmission rate in \eqref{eq12}. $ \hfill \blacksquare $

\begin{lemma}  \label{lemma2}
	The CDF and PDF of the transmission delay in the single-antenna scenario are given by
	\begin{equation*} 
		F_{\tau}(t) = e^{-\frac{1}{\rho} \left( 2^{\frac{L}{Bt}} - 1 \right)}, \tag{14} \label{eq14}
	\end{equation*}
	\begin{equation*} 
		f_{\tau}(t) = \frac{L \ln 2}{\rho B} \frac{1}{t^2} 2^{\frac{L}{Bt}} e^{-\frac{1}{\rho} \left( 2^{\frac{L}{Bt}} - 1 \right)}. \tag{15} \label{eq15}
	\end{equation*}
\end{lemma}

\textit{Proof:}
\begin{align*} 
	F_{\tau}(t) &= \Pr \left\{ \frac{L}{BR} \le t \right\}   \\
	& = \Pr \left\{ R \ge \frac{L}{Bt} \right\} \\
	& = 1 - F_R \left( \frac{L}{Bt} \right). \tag{16} \label{eq16}
\end{align*}
Substituting \eqref{eq11} into the above equation yields the CDF of the transmission delay in \eqref{eq14}. The derivation of the CDF yields the PDF of the transmission delay in \eqref{eq15}.$ \hfill \blacksquare $

\subsection{Delay Distribution in the Multi-Antenna Scenario} \label{sec3.2}
Multi-antenna transmission can effectively exploit spatial degrees of freedom to improve transmission delay performance. To illustrate the impact of the number of antennas on delay performance, we derive the CDF and PDF of transmission rate and transmission delay based on \eqref{eq5}, which are given by Lemma \ref{lemma3} and Lemma \ref{lemma4}, respectively.

\begin{lemma}  \label{lemma3}
	The CDF and PDF of the transmission rate in the multi-antenna scenario are given by
	\begin{equation*} 
		F_{R}(y) = \frac{\gamma \left(N,\frac{N}{\rho}\left(2^y-1\right)\right)}{\Gamma(N)}, \tag{17} \label{eq17}
	\end{equation*}
	\begin{equation*} 
		f_{R}(y) = \frac{N^N \ln 2}{\rho^N \Gamma(N)} 2^y \left(2^y-1\right)^{N-1} e^{-\frac{N}{\rho} \left( 2^y - 1 \right) }. \tag{18} \label{eq18}
	\end{equation*}
\end{lemma}

\textit{Proof:}
The proof of Lemma \ref{lemma3} is similar to that of Lemma \ref{lemma1} and is omitted here. $ \hfill \blacksquare $

\begin{lemma}  \label{lemma4}
	The CDF and PDF of the transmission delay in the multi-antenna scenario can be represented as
	\begin{equation*} 
		F_{\tau}(t) = \frac{\Gamma\left(N,\frac{N}{\rho}\left(2^{\frac{L}{Bt}} - 1\right)\right)}{\Gamma(N)} , \tag{19} \label{eq19}
	\end{equation*}
	\begin{equation*} 
		f_{\tau}(t) = \frac{ LN^N\ln 2 }{B\rho^N \Gamma(N)} \frac{1}{t^2} 2^{\frac{L}{Bt}} \left( 2^{\frac{L}{Bt}} - 1 \right)^{N-1} e^{-\frac{N}{\rho}\left( 2^{\frac{L}{Bt}} - 1 \right)} , 
		\tag{20} \label{eq20} 
	\end{equation*}   
	where $ \Gamma(\cdot, \cdot) $ is the upper incomplete gamma function \cite[eq. (8.350.2)]{Gradshteyn} and $ \Gamma(s, x) + \gamma(s, x) = \Gamma(s) $.
\end{lemma}

\textit{Proof:}
The proof of Lemma \ref{lemma4} is similar to that of Lemma \ref{lemma2} and is omitted here. $ \hfill \blacksquare $
 
\begin{remark}  \label{remark1}
	It can be observed that single-antenna systems are a special case of multi-antenna systems, and Lemma \ref{lemma1} and Lemma \ref{lemma2} can be obtained by setting $N = 1$ in Lemma \ref{lemma3} and Lemma \ref{lemma4}. The conclusion is also applicable to finite blocklength transmission. Therefore, no further distinction is made in the subsequent analysis, and multi-antenna transmission is used as an example for a comprehensive analysis.
\end{remark}

\begin{remark}  \label{remark2}
	The delay violation probability is defined as the probability that the delay exceeds the threshold value $\tau_{th}$, i.e., $p_{v} \left( \tau_{th} \right) = \Pr \left\{ \tau > \tau_{th} \right\} $. According to Lemma \ref{lemma4}, we have
	\begin{equation*} 
		p_{v} \left( \tau_{th} \right) = 1 - F_{\tau} \left( \tau_{th} \right) = \frac{\gamma\left(N,\frac{N}{\rho}\left(2^{\frac{L}{B\tau_{th}}} - 1\right)\right)}{\Gamma(N)}. \tag{21} \label{eq21}
	\end{equation*}
\end{remark}

\subsection{Average Delay and Jitter} \label{sec3.3}

In the following discussion, we will focus on jitter. According to \eqref{eq10}, the moment of the delay can be formulated as 
\begin{equation}
	\small
	\mathbb{E} [\tau^i] = \frac{LN^N\ln 2}{B \rho^N \Gamma(N)} \int_0^{\infty} t^{i-2} 2^{\frac{L}{Bt}} \left(2^{\frac{L}{Bt}}-1\right)^{N-1} e^{-\frac{N}{\rho}\left( 2^{\frac{L}{Bt}}-1 \right)} dt. 
	\tag{22} \label{eq22} 
\end{equation}
	\normalsize
It can be seen that solving \eqref{eq22} is very difficult due to the coupling relationship between exponential and power functions, which makes it challenging to obtain a closed-form solution. Note that, according to the scaling properties of expectation and variance, the average delay and jitter can also be expressed as
\begin{equation*} 
	\mathbb{E} [\tau ] = \frac{L}{B} \mathbb{E} [R^{-1}], \tag{23}  \label{eq23}
\end{equation*}
\begin{equation*} 
	\text{Var} [\tau ] = \frac{L^2}{B^2} \text{Var} [R^{-1}]. \tag{24} \label{eq24}
\end{equation*}
Therefore, we analyze the average delay and jitter based on \eqref{eq23} and \eqref{eq24}, and the analytical results are shown in Theorem \ref{theorem1} and Theorem \ref{theorem2}.

\begin{theorem}  \label{theorem1}
	The average delay and jitter can be approximated as
	\begin{equation*}
		\mathbb{E} [\tau ] \approx \frac{L}{B} \frac{\mathbb{E}[R^2]}{\mathbb{E}[R]^3}, \tag{25} \label{eq25}
	\end{equation*}
	\begin{equation*} 
		\text{Var} [\tau ] \approx \frac{L^2}{B^2} \left[ -\frac{\mathbb{E}[R^2]^2}{\mathbb{E}[R]^6} + \frac{3\mathbb{E}[R^2]}{\mathbb{E}[R]^4} - \frac{2}{\mathbb{E}[R]^2} \right], \tag{26} \label{eq26}
	\end{equation*}
	where $\mathbb{E}\left[R\right]$ and $\mathbb{E}\left[R^2\right]$ are the first-order and second-order moments of the transmission rate, respectively. $\mathbb{E}\left[R\right]$ and $\mathbb{E}\left[R^2\right]$ can be computed under the Rayleigh channel by
	\begin{align*} 
		\mathbb{E} [R] = &
		\frac{N^N e^{\frac{N}{\rho}}}{\ln2 \rho^{N}\Gamma(N)} \sum_{k=0}^{N-1} a_kb \Bigg[ \Psi\left(\frac{1}{b}+k+1,\frac{N}{\rho}\right) \\&
		- \Psi\left(k+1,\frac{N}{\rho}\right) \Bigg], \tag{27} \label{eq27}
	\end{align*}
	\begin{align*} 
		\mathbb{E} [R^2] =& \frac{N^N e^{\frac{N}{\rho}}}{(\ln2)^2\rho^{N}\Gamma(N)} \sum_{k=0}^{N-1} a_k b^2  \Bigg[\Psi\left(\frac{2}{b}+k+1,\frac{N}{\rho}\right) \\
		& - 2\Psi\left(\frac{1}{b}+k+1,\frac{N}{\rho}\right) + \Psi\left(k+1,\frac{N}{\rho}\right) \Bigg], \tag{28} \label{eq28}
	\end{align*}
	where $ a_k = \begin{pmatrix} N-1\\k \end{pmatrix} (-1)^{N-1-k} $, $ b $ is a large constant, and $ \Psi(x,y) \triangleq y^{-x}\Gamma(x,y) $.
\end{theorem}

\textit{Proof:}
Please refer to Appendix \ref{appendix1} for details. $ \hfill \blacksquare $

\begin{theorem}  \label{theorem2}
	The average delay and jitter can also be approximated as
	\begin{equation*} 
		\mathbb{E} [\tau] \approx \frac{L}{B} \left[ \frac{1}{ \mathbb{E} [R] } + \frac{\text{Var}[R]}{(\mathbb{E}[R])^3} \right], \tag{29} \label{eq29} 
	\end{equation*}
	\begin{equation*} 
		\text{Var}  [\tau] \approx \frac{L^2}{B^2} \left[ \frac{\text{Var}[R]}{ (\mathbb{E}[R])^4 } - \frac{ (\text{Var}[R])^2 }{(\mathbb{E}[R])^6} \right], \tag{30} \label{eq30}
	\end{equation*}
	where $ \text{Var}\left[R\right] $ is the variance of the transmission rate. $\mathbb{E}\left[R\right]$ and $ \text{Var}\left[R\right] $ can be approximated under the Rayleigh channel by the following equations: 
	\begin{equation*}  
		\mathbb{E} [R] \approx \log_2 \left(1+\rho\right) 
		- \frac{\rho^2}{2\ln2N \left(1+ \rho\right)^2}, \tag{31} \label{eq31} 
	\end{equation*}
	\begin{equation*} 
		\text{Var} [R] \approx \frac{\rho^2}{(\ln2)^2 N \left(1+\rho\right)^2} 
		- \frac{\rho^4}{4(\ln2)^2 N^2 \left(1+ \rho\right)^4}. \tag{32} \label{eq32}
	\end{equation*}
\end{theorem}

\textit{Proof:}
We can think of the transmission rate $R=\log_2(1+\gamma)$ as a function of the random variable $\gamma$. According to Appendix \ref{appendix1} and the mean and variance of the gamma random variable, we can derive \eqref{eq31} and \eqref{eq32}. $ \hfill \blacksquare $

\section{Finite Blocklength Analysis} \label{sec4}
For the sake of simplicity, we did not consider retransmission.\footnote[3]{On the one hand, based on adaptive coding and modulation (ACM) technology, the transmitter is able to choose an appropriate blocklength during each transmission process so that the transmission rate does not exceed the channel capacity, thus achieving high reliability. On the other hand, the retransmission delays introduced by different retransmission mechanisms are different and need to be studied separately. Therefore, we leave the retransmission delay as future work.} We assume that the block error rate (BLER) is $ \epsilon $, and the maximum achievable coding rate can be inscribed by the approximate formula proposed by Polyanskiy et al. \cite{Polyanskiy}, which is  
\begin{equation*} 
	R (\gamma, n, \epsilon) \approx \log_2 \left( 1+\gamma \right)  - \sqrt{\frac{V(\gamma)}{n}} Q^{-1}(\epsilon), \tag{33} \label{eq33}
\end{equation*}
where $\gamma$ represents the instantaneous SNR, $V(\gamma) = \left[ 1- \frac{1} {( 1+\gamma )^2}\right] (\log_2 e)^2$ denotes the channel dispersion, and $ Q^{-1}\left( x \right) $ is the inverse function of $ Q \left( x \right) = \int_{x}^{\infty} \frac{1}{\sqrt{2\pi}} e^{-\frac{t^2}{2}}\, {\rm d}t $. The number of information bits that can be transmitted should satisfy
\begin{equation*}  
	n\log_2 \left( 1+\gamma\right) - \sqrt{nV(\gamma)} Q^{-1}\left( \epsilon \right) \ge L. \tag{34} \label{eq34}
\end{equation*}

We define 
\begin{equation*} 
	F(n,\gamma) \triangleq n\log_2 \left( 1+\gamma\right) - \sqrt{nV(\gamma)} Q^{-1}\left( \epsilon \right) - L. \tag{35} \label{eq35}
\end{equation*}
It is clear that the minimum blocklength required for transmitting $L$ bits is an implicit function determined by $F(n,\gamma) = 0$. Solving this quadratic equation yields  
\begin{align*} 
	\small
	&n(\gamma,L,\epsilon) =\\
	& \left(\frac{\sqrt{V(\gamma)} Q^{-1}\left( \epsilon \right) + \sqrt{V(\gamma) \left(Q^{-1}\left( \epsilon \right)\right)^2 + 4L\log_2 \left( 1+\gamma\right)} }{2\log_2 \left( 1+\gamma\right)} \right)^2. \tag{36} \label{eq36}
\end{align*}
\normalsize
The transmission delay is expressed as $\tau = \frac{n}{B}$ \cite{Durisi}.

\subsection{Delay Distribution} \label{sec4.1}
We derive the CDF and PDF of the transmission delay based on the CDF expression in \eqref{eq5} and the coding rate in \eqref{eq33}, as given by Lemma \ref{lemma5}.

\begin{lemma}  \label{lemma5} 
	The CDF and PDF of the transmission delay are given by
	\begin{equation*}  
		F_{\tau}(t) = \frac{\Gamma \left(N,\frac{N}{\rho}\left(e^{u(t)}-1\right) \right)}{\Gamma(N)}, 
		\tag{37} \label{eq37} 
	\end{equation*}
	\begin{equation*}  
		f_{\tau}(t) = -\frac{N^N}{\rho^N \Gamma(N)} u'(t) e^{u(t)} \left[ e^{u(t)} - 1 \right]^{N-1} e^{-\frac{N}{\rho}\left( e^{u(t)} - 1 \right)} , \tag{38} \label{eq38} 
	\end{equation*}
	where
	\begin{align*}  
		u(t) =& \frac{Q^{-1}\left( \epsilon \right)}{\sqrt{Bt}} + \frac{L\ln2}{Bt} \\
		&- \frac{1}{2} \sum_{m=1}^M \frac{m^{m-1}}{m!}\left( \frac{Q^{-1}\left(\epsilon\right) }	{ \sqrt{Bt} } e^{ - \frac{2Q^{-1}\left( \epsilon \right)} {\sqrt{Bt}} - \frac{2L\ln2}{Bt}  } \right)^{m}, \tag{39} \label{eq39}
	\end{align*}
	\begin{small}
		\begin{align*}  
			&u'(t) = - \frac{Q^{-1}\left( \epsilon \right)}{2\sqrt{Bt^3}} - \frac{L\ln2}{Bt^2}  - \frac{1}{2} \sum_{m=1}^M \frac{m^{m-1}}{(m-1)!} \\
			&\times \left( \frac{Q^{-1}\left(\epsilon\right) }	{ \sqrt{Bt} } e^{ - \frac{2Q^{-1}\left( \epsilon \right)} {\sqrt{Bt}} - \frac{2L\ln2}{Bt} } \right)^{m} \left( \frac{Q^{-1}\left( \epsilon \right)}{\sqrt{Bt^3}} + \frac{2L\ln2}{Bt^2} - \frac{1}{2t} \right). \tag{40} \label{eq40}
		\end{align*}
	\end{small}It is noted that the infinite series summing in the third terms of $u(t)$ and $u'(t)$ have high computational complexity. We approximate the third term by summing the first $M$ terms. The larger the value of $M$, the more accurate the approximation.
\end{lemma}

\textit{Proof:}
Please refer to Appendix \ref{appendix2} for details. $ \hfill \blacksquare $

\begin{remark}  \label{remark3}
	Let $N=1$, and the CDF and PDF of transmission delay for the single-antenna scenario can be obtained as
	\begin{equation*}  
		F_{\tau}(t) = e^{-\frac{1}{\rho} \left( e^{u(t)} - 1 \right) }, \tag{41} \label{eq41}
	\end{equation*}
	\begin{equation*}  
		f_{\tau}(t) = -\frac{1}{\rho} u'(t) e^{u(t)}  e^{-\frac{1}{\rho} \left( e^{u(t)} - 1 \right) }. \tag{42} \label{eq42}
	\end{equation*}
\end{remark}

Due to the complexity of the PDF expression in \eqref{eq38}, it is impossible to derive the moments of the delay similar to infinite blocklength transmission. Therefore, we attempt to obtain an approximation expression for the delay at high SNR.

\subsection{Approximation Expression of Delay for High SNR}  \label{sec4.2}
At high SNR, the channel dispersion can be approximated as $ V(\gamma) \approx (\log_2 e)^2 $ \cite{Feng}, and thus, the transmission delay can be expressed as
\begin{equation*}  \tau \approx \left(\frac{ Q^{-1}\left( \epsilon \right) + \sqrt{ \left(Q^{-1}\left( \epsilon \right)\right)^2 + 4 L\ln 2 \ln \left( 1+\gamma\right) } }{2\sqrt{B}\ln \left( 1+\gamma\right)} \right)^2 . \tag{43} \label{eq43}
\end{equation*}
The above approximation is sufficiently accurate when SNR is greater than 5 dB, and this SNR condition is usually easy to meet in general wireless communication systems. This approximation is widely used in research related to URLLC \cite{Feng}, \cite{C. Sun}. Moreover, the approximation provides an upper bound on the transmission delay, which imposes stricter requirements on the transmission delay. We obtain the approximate expressions for the CDF and PDF of the transmission delay in Lemma \ref{lemma6}.

\begin{lemma}  \label{lemma6}
	The approximate expressions for the CDF and PDF of the transmission delay at high SNR are given by 
	\begin{equation*}  F_{\tau}(t) = \frac{\Gamma \left(N,\frac{N}{\rho}\left( e^{\frac{Q^{-1}\left(\epsilon \right)}{\sqrt{Bt}} + \frac{L\ln2}{Bt} }-1 \right) \right)}{\Gamma(N)} , \tag{44} \label{eq44}
	\end{equation*}
	\begin{align*}  f_{\tau}(t) =& \frac{N^N}{\rho^N \Gamma(N)}  \left( \frac{Q^{-1}\left(\epsilon\right)}{2\sqrt{Bt^3}} + \frac{L\ln2}{Bt^2} \right)
		e^{\frac{Q^{-1}\left(\epsilon \right)}{\sqrt{Bt}} + \frac{L\ln 2}{Bt} }  \\
		& \times \left( e^{\frac{Q^{-1}\left(\epsilon \right)}{\sqrt{Bt}} + \frac{L\ln 2}{Bt} }-1 \right)^{N-1} e^{-\frac{N}{\rho}\left( e^{\frac{Q^{-1}\left(\epsilon \right)}{\sqrt{Bt}} + \frac{L\ln2}{Bt} }-1 \right)}. \tag{45} \label{eq45}
	\end{align*}
\end{lemma}

\textit{Proof:}
The proof is similar to that of Lemma \ref{lemma5} and is omitted here. $ \hfill \blacksquare $

\begin{remark}  \label{remark4}
	Similarly to infinite blocklength transmission, the delay violation probability for finite blocklength transmission can be obtained as
	\begin{equation*} 
		p_{v} \left( \tau_{th} \right) = \frac{\gamma \left(N,\frac{N}{\rho}\left( e^{\frac{Q^{-1}\left(\epsilon \right)}{\sqrt{B\tau_{th} }} + \frac{L\ln2}{B\tau_{th} } }-1 \right) \right)}{\Gamma(N)} . \tag{46} \label{eq46}
	\end{equation*}
\end{remark}

The expression in \eqref{eq45} is still too complex for theoretical analysis, and it is nearly impossible to derive the average delay and jitter from it. Therefore, we simplify and approximate \eqref{eq43} and then derive a tight upper bound of delay in Lemma \ref{lemma7}.

\begin{lemma}  \label{lemma7}
	The approximation and upper bound of the transmission delay at high SNR are given by
	\begin{equation*}  \tau \approx \frac{ L\ln 2 }{ B\ln \left( 1+\gamma\right) }
		+ \frac{\sqrt{L\ln2} Q^{-1}\left( \epsilon \right)}{B \left[\ln(1+\gamma)\right]^{\frac{3}{2}}  } 
		+ \frac{\left(Q^{-1}\left( \epsilon \right) \right)^2}{ 2B\left[\ln\left( 1+\gamma\right)\right]^2 }, \tag{47} \label{eq47}
	\end{equation*}
	\begin{equation*}  \tau^{\text{upper}} = \frac{ 2L\ln 2+\left(Q^{-1}\left( \epsilon \right) \right)^2 }{ 2B\ln \left( 1+\gamma\right) }
		+ \frac{\sqrt{L\ln2} Q^{-1}\left( \epsilon \right)}{B \left[\ln(1+\gamma)\right]^{\frac{3}{2}} }. \tag{48} \label{eq48}
	\end{equation*}
\end{lemma}

\textit{Proof:}
\eqref{eq43} can be transformed into
\begin{align*}  \tau =& \frac{ L\ln 2 }{ B\ln \left( 1+\gamma\right) }
	+ \frac{\left(Q^{-1}\left( \epsilon \right) \right)^2}{ 2B\left[\ln \left( 1+\gamma\right)\right]^2 } \\&+ \frac{\sqrt{L\ln2} Q^{-1}\left( \epsilon \right)}{B \left[\ln(1+\gamma)\right]^{\frac{3}{2}} } \sqrt{ 1 + \frac{\left(Q^{-1}\left( \epsilon \right)\right)^2}{4L\ln2 \ln\left( 1+\gamma\right) } }. \tag{49} \label{eq49}
\end{align*}
When $\gamma \ge e^{\frac{\left(Q^{-1}\left(\epsilon\right)\right)^2}{4L\ln2} } - 1$, we have $4L\ln2 \ln\left(1+\gamma\right) \ge \left(Q^{-1}\left(\epsilon\right)\right)^2 \ge 0 $. Since the inverse function of the Gaussian Q-function is a monotonically decreasing function, when $ \epsilon = 10^{-9}$, $L = 100$, $e^{\frac{\left(Q^{-1}\left(\epsilon\right)\right)^2}{4L\ln2} } - 1 \approx 0.14$, $\gamma \ge e^{\frac{\left(Q^{-1}\left(\epsilon\right)\right)^2}{4L\ln2} } - 1$ clearly holds for the high SNR case considered in this subsection, and therefore, $4L\ln2 \ln\left(1+\gamma\right) \ge \left(Q^{-1}\left(\epsilon\right)\right)^2 $ also holds. The binomial approximation of $\sqrt{ 1 + \frac{\left(Q^{-1}\left( \epsilon \right)\right)^2}{4L\ln2 \ln\left( 1+\gamma\right) } }$ yields
\begin{align*}  \tau \approx& \frac{ L\ln 2 }{ B\ln \left( 1+\gamma\right) }
	+ \frac{\sqrt{L\ln2} Q^{-1}\left( \epsilon \right)}{B\left[\ln(1+\gamma)\right]^{\frac{3}{2}}  } 
	+ \frac{\left(Q^{-1}\left( \epsilon \right) \right)^2}{ 2B\left[\ln \left( 1+\gamma\right)\right]^2 } \\ &	+ \frac{\left(Q^{-1}\left( \epsilon \right)\right)^3 }{8B \sqrt{L\ln2} \left[\ln(1+\gamma)\right]^{\frac{5}{2}}}. \tag{50} \label{eq50}
\end{align*}
The fourth term in \eqref{eq50} can be ignored because it is smaller than the first three terms. At this point, the approximation in \eqref{eq47} for the transmission delay at high SNR can be obtained. Furthermore, under the high SNR condition, since $\ln(1+\gamma) > 1$, we have $\frac{1}{\left[\ln(1+\gamma)\right]^2} < \frac{1}{\left[\ln(1+\gamma)\right]^{\frac{3}{2}}} < \frac{1}{\ln(1+\gamma)}$. Additionally, since the third term in \eqref{eq50} is slightly smaller than the first two, the upper bound on the delay at high SNR can be obtained from \eqref{eq48}. $ \hfill \blacksquare $

\subsection{Average Delay and Jitter for High SNR}  \label{sec4.3}
\begin{theorem}  \label{theorem3}
	At high SNR, we can approximate the average delay and jitter as
	\begin{align*}  \mathbb{E}\left[ \tau \right] \approx &
		\frac{ 2L\ln 2+\left(Q^{-1}\left( \epsilon \right) \right)^2 }{ 2\ln2B } \left[ \frac{1}{ \mathbb{E} [R] } + \frac{\text{Var}[R]}{(\mathbb{E}[R])^3} \right]  \\
		& + \frac{\sqrt{L\ln2} Q^{-1}\left( \epsilon \right)}{B(\ln2)^{\frac{3}{2}}} \left[ \frac{1}{ (\mathbb{E}[R])^{\frac{3}{2}} } + \frac{ 15\text{Var}[R]}{8(\mathbb{E}[R])^{\frac{7}{2}} } \right], \tag{51} \label{eq51}
	\end{align*}
	\begin{small}
		\begin{align*}  \text{Var} \left[ \tau \right] \approx &
			\left[\frac{ 2L\ln 2+\left(Q^{-1}\left( \epsilon \right) \right)^2 }{ 2\ln2 B}\right]^2  
			\left[ \frac{\text{Var}[R]}{ (\mathbb{E}[R])^4 } - \frac{ (\text{Var}[R])^2 }{(\mathbb{E}[R])^6} \right] \\
			& + \frac{L\ln2 \left(Q^{-1}\left( \epsilon \right)\right)^2}{B^2\left(\ln2\right)^3 } \left[\frac{9 \text{Var} [R] }{4\left(\mathbb{E}[R]\right)^5} 
			- \frac{225 \left(\text{Var} [R]\right)^2 }{64\left(\mathbb{E}[R]\right)^7} \right]
			, \tag{52} \label{eq52}
		\end{align*}
	\end{small}where $ \mathbb{E}\left[R\right] $ and $ \text{Var}\left[R\right] $ can be obtained from either Theorem \ref{theorem1} or Theorem \ref{theorem2} under the Rayleigh channel.
\end{theorem}

\textit{Proof:}
Please refer to Appendix \ref{appendix3} for details. $ \hfill \blacksquare $

\section{Simulation Results and Discussions} \label{sec5}
In this section, we validate the accuracy of the theoretical analysis and approximation results through simulations and numerical calculations. Unless otherwise stated, the simulation parameters are listed in Table \ref{table1}.

\begin{table} 
	\caption{SIMULATION PARAMETERS}
	\begin{tabular}{ll}
		\hline
		Parameter                                 & Defult value                                                                    \\ \hline
		Average SNR ($\rho$)                      & 10 dB                                                                            \\
		The number of information bits ($L$)      & \begin{tabular}[c]{@{}l@{}}1000 bits for IBL and \\ 200 bits for FBL\end{tabular} \\
		The number of channel realizations        & $10^6$                                                                          \\
		The number of transmitter antennas ($N$)  & 8                                                                               \\
		Bandwidth ($B$)                           & 200 kHz                                                                          \\
		BLER ($\epsilon$)                         & $10^{-7}$                                                                \\ \hline
	\end{tabular}
	\label{table1}
\end{table}

\subsection{Accuracy of the Approximation in \eqref{eq58} and \eqref{eq64}} \label{sec5.1}
\begin{figure}[htbp]    
	\centering       
	\includegraphics[width=0.43\textwidth]{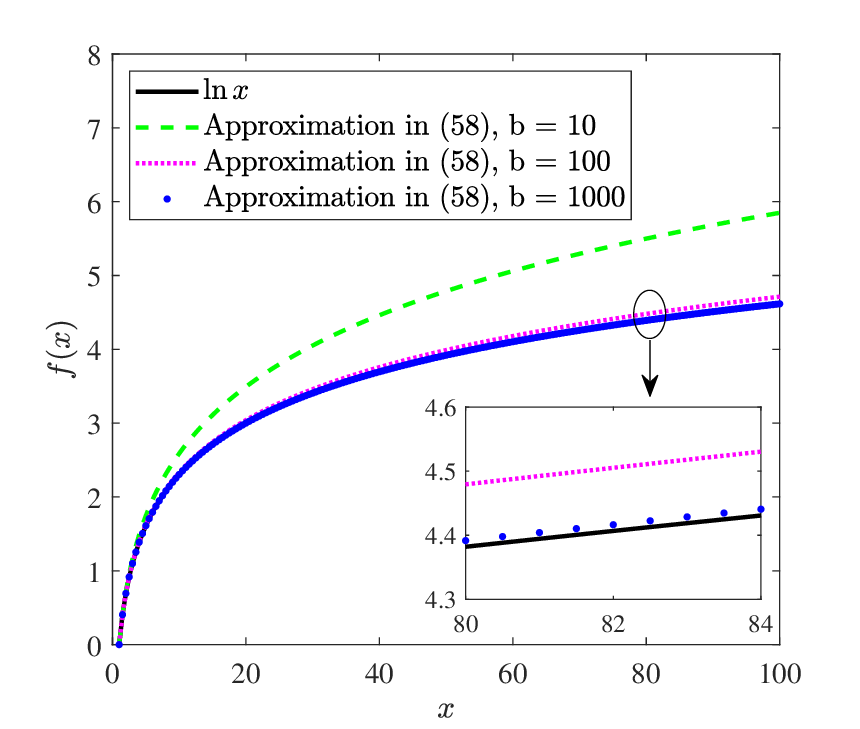} 
	\caption{The tightness of the approximation in (\ref{eq58}) in Appendix \ref{appendix1}.}        
	\label{fig.1}                      
\end{figure}

\begin{figure}[htbp]    
	\centering       
	\includegraphics[width=0.43\textwidth]{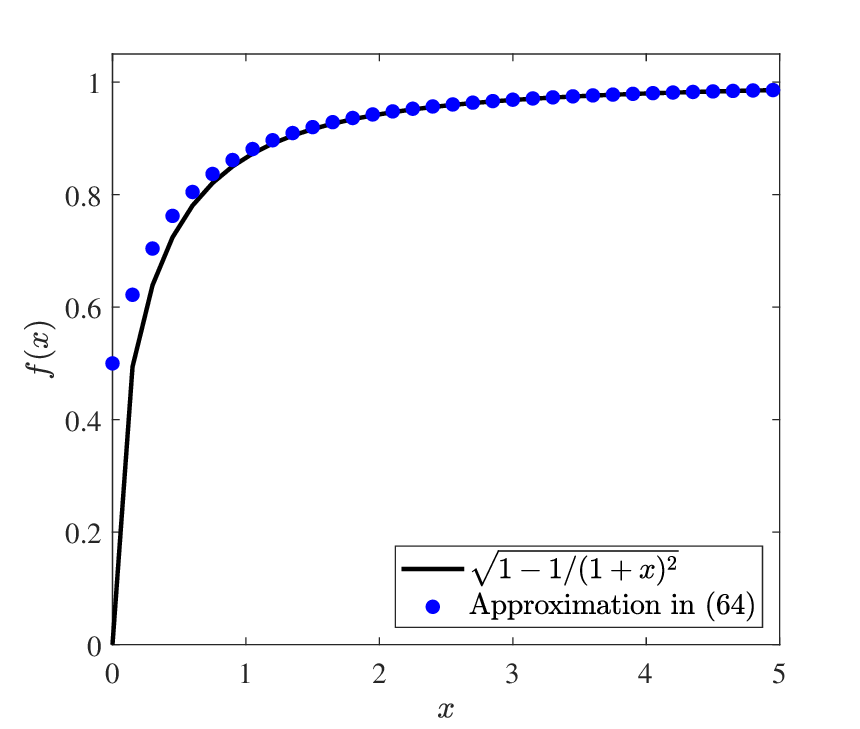} 
	\caption{The tightness of the approximation in (\ref{eq64}) in Appendix \ref{appendix2}.}      
	\label{fig.2}                      
\end{figure}

Figures \ref{fig.1} and \ref{fig.2} illustrate the tightness of the approximations in \eqref{eq58} and \eqref{eq64}, respectively. As shown in Fig. \ref{fig.1}, as the constant $b$ increases, the approximation in \eqref{eq58} perfectly matches with $\ln x$, which indicates that the approximation is accurate when a large value of $b$ is chosen, such as $b = 1000$. As a result, the calculations for $\mathbb{E}[R]$ and $\text{Var}[R]$ in Appendix \ref{appendix1} are accurate. Figure \ref{fig.2} shows that the curve of eq. \eqref{eq64} is close to the curve of $ \sqrt{1-1/(1+x)^2} $ when $x > 1$, indicating that the approximation in Appendix \ref{appendix2} is highly accurate when SNR is large. This suggests that the computation of the CDF and PDF of the transmission delay in the FBL regime is accurate.  

\subsection{CDF and PDF of Transmission Rate and Transmission Delay} \label{sec5.2}
\begin{figure}[htbp] 
	\centering  
	\subfigure[]{
		\label{fig3.sub.1}
		\includegraphics[width=0.43\textwidth]{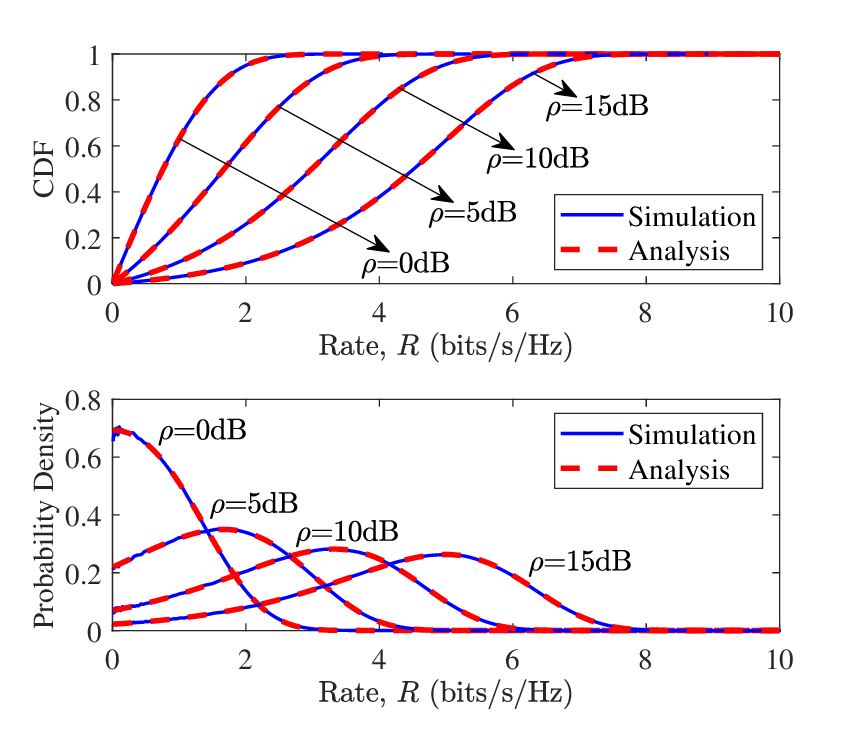}}	
	\subfigure[]{
		\label{fig3.sub.2}
		\includegraphics[width=0.43\textwidth]{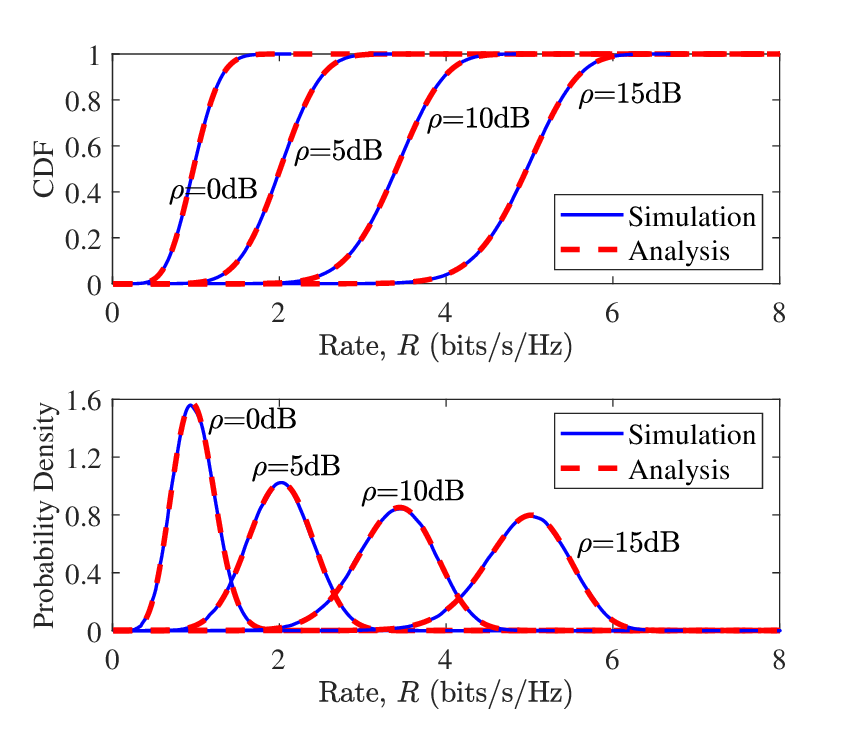}}
	\caption{CDF and PDF of transmission rate under infinite blocklength regime with $L = 1000$ bits. (a) Single-antenna scenario. (b) Multi-antenna scenario.}
	\label{fig.3}
\end{figure}

In Fig. \ref{fig.3}, we compare the theoretical and simulation results of the transmission rate in the IBL regime. The theoretical results agree with their corresponding simulation results, verifying the correctness of the transmission rate distribution expressions in the IBL regime (i.e., Lemma \ref{lemma1} and Lemma \ref{lemma3}). Figure \ref{fig3.sub.1} demonstrates the CDF and PDF of the transmission rate for the single-antenna scenario. It can be seen that the transmission rate has a large probability density at 0 bits/s/Hz when the average SNR $\rho$ is low. Furthermore, as the average SNR increases, the PDF curves gradually shift to the right, and the probability density at 0 bits/s/Hz gradually decreases. Figure \ref{fig3.sub.2} shows the CDF and PDF of the transmission rate in the multi-antenna scenario. The probability density of the transmission rate at 0 bits/s/Hz is close to zero. The main reason is that the channels between different antennas are independent, which reduces the probability of deep fading.

\begin{figure}[htbp]
	\centering  
	\subfigure[]{
		\label{fig4.sub.1}
		\includegraphics[width=0.43\textwidth]{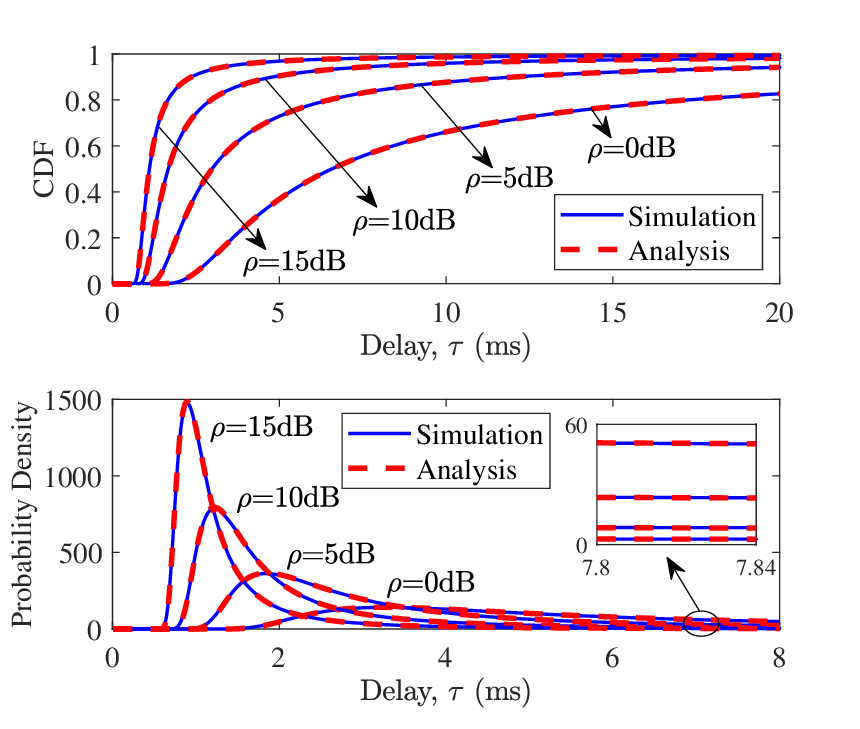}}	
	\subfigure[]{
		\label{fig4.sub.2}
		\includegraphics[width=0.43\textwidth]{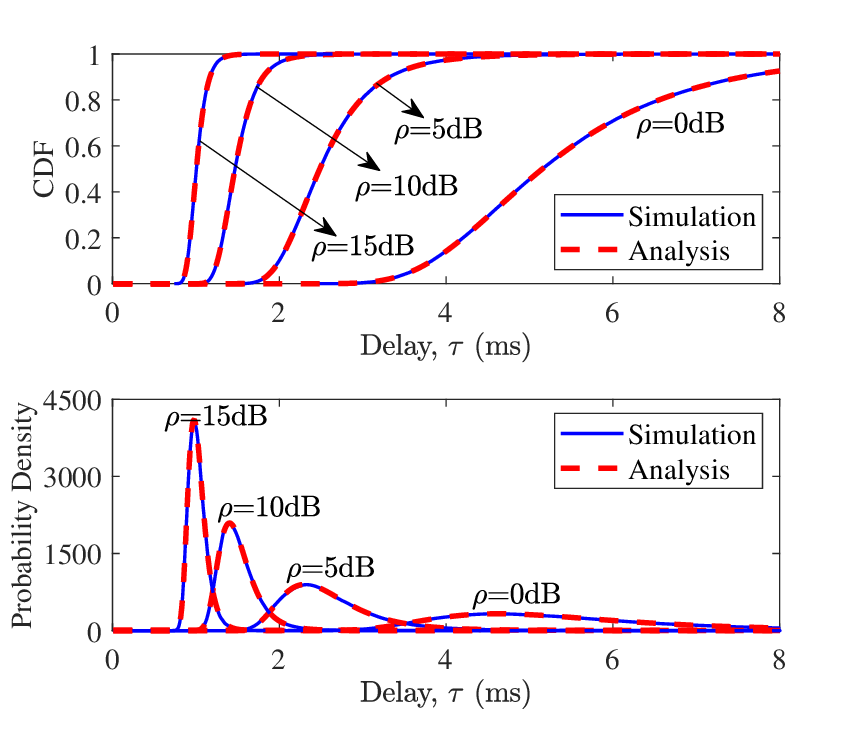}}
	\caption{CDF and PDF of transmission delay under infinite blocklength regime with $L = 1000$ bits. (a) Single-antenna scenario. (b) Multi-antenna scenario.}
	\label{fig.4}
\end{figure}

Figure \ref{fig.4} presents the theoretical and simulation results of the transmission delay in the IBL regime. We can see that the theoretical and simulation results are in good agreement, which verifies the correctness of the transmission delay distribution expressions in the IBL regime (i.e., Lemma \ref{lemma2} and Lemma \ref{lemma4}). Figure \ref{fig4.sub.1} presents the CDF and PDF of the transmission delay for the single-antenna scenario. It is worth noting that the PDF curves have a large tail distribution. This is because when the channel experiences deep fading, the channel capacity becomes close to zero due to limited transmit power and bandwidth, resulting in a large delay. Furthermore, a smaller average SNR will result in a more significant tail distribution since the probability density of the transmission rate at 0 bits/s/Hz is higher at low SNR region. As the average SNR increases, the probability density of the tail distribution decreases because the probability density of the transmission rate at 0 bits/s/Hz decreases. From these results, we see that the transmission delay distribution in the single-antenna scenario is not satisfactory at low SNR. Furthermore, because of the tail distribution of transmission delay, the average delay and jitter will be infinite. Figure \ref{fig4.sub.2} illustrates the CDF and PDF of the transmission delay for the multi-antenna scenario. The tail distribution is significantly improved in the multi-antenna scenario. This is reasonable because multi-antenna transmission reduces the probability of deep fading, resulting in a decrease in the probability density of transmission rate at 0 bits/s/Hz. However, the delay distribution remains unfavorable at low SNR, e.g., $\rho = 0$ dB. Compared with Fig. \ref{fig4.sub.1}, the PDF curves in Fig. \ref{fig4.sub.2} become narrower. This indicates that the multi-antenna transmission is beneficial for reducing jitter. The PDF curves become narrower as the average SNR increases, which indicates that increasing the transmit power reduces jitter.


\begin{figure}[htbp]  
	\centering  
	\subfigure[]{
		\label{fig5.sub.1}
		\includegraphics[width=0.43\textwidth]{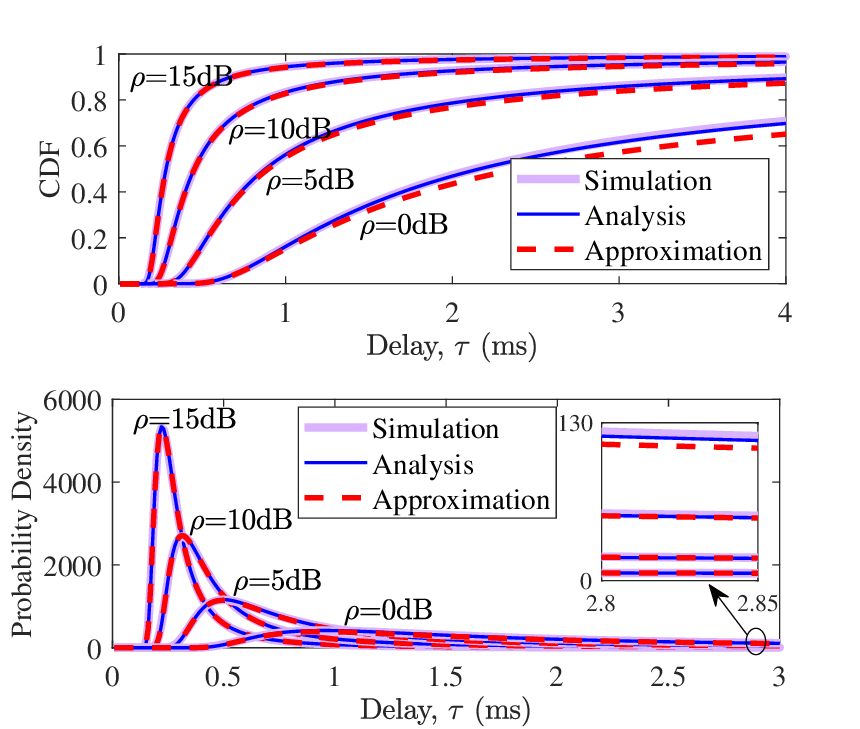}}	
	\subfigure[]{
		\label{fig5.sub.2}
		\includegraphics[width=0.43\textwidth]{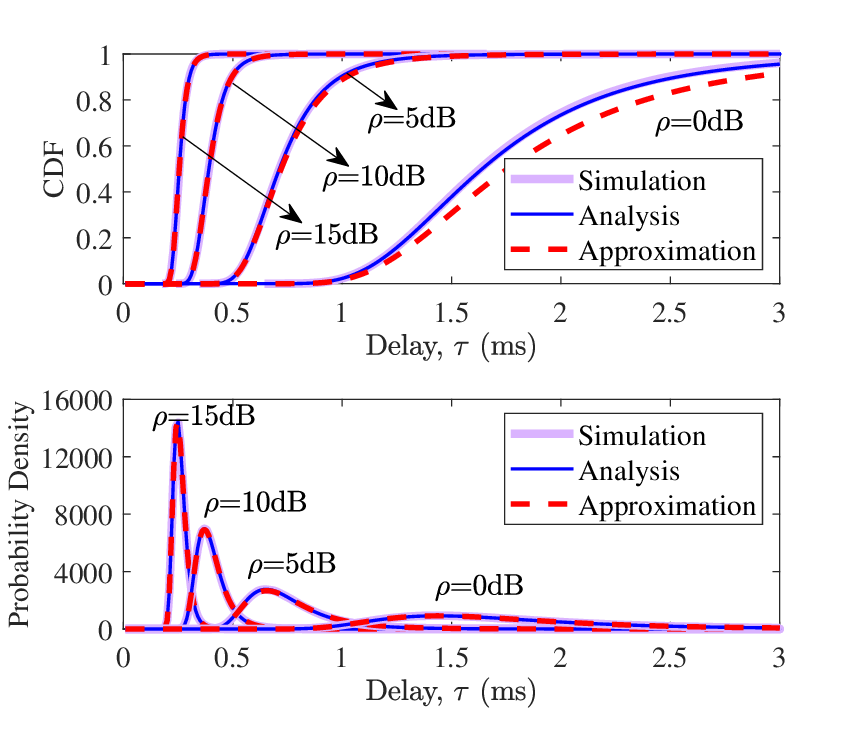}}
	\caption{CDF and PDF of transmission delay under finite blocklength regime with $L = 200 $ bits. (a) Single-antenna scenario. (b) Multi-antenna scenario.}
	\label{fig.5}
\end{figure}

In Fig. \ref{fig.5}, we show the theoretical, high SNR approximation, and simulation results of transmission delay in the FBL regime. The agreement among the three curves proves the correctness of the transmission delay distribution expressions in the FBL regime (i.e., Lemma \ref{lemma5} and Lemma \ref{lemma6}). Figure \ref{fig5.sub.1} depicts the CDF and PDF of the transmission delay for the single-antenna scenario. Comparison of Fig. \ref{fig5.sub.1} with Fig. \ref{fig4.sub.1} indicates that the transmission delay is larger in the FBL case in terms of transmitting the same number of information bits. This is reasonable since the maximum coding rate in the FBL regime is lower than Shannon capacity. Therefore, the delay performance will be overestimated if the Shannon capacity is considered. Notice that the tail distribution still exists and is more severe than in the IBL case. The CDF and PDF of the transmission delay for the multi-antenna scenario are shown in Fig. \ref{fig5.sub.2}. Comparing Fig. \ref{fig5.sub.1} and Fig. \ref{fig5.sub.2}, it can be observed that the transmission delay distribution is more desirable in the multi-antenna scenario, which suggests that multi-antenna transmit diversity can enhance delay performance. These results are consistent with the conclusions obtained in Fig. \ref{fig.4}.

\subsection{Average Delay and Jitter} \label{sec5.3}

\begin{figure}[htbp]     
	\centering       
	\includegraphics[width=0.43\textwidth]{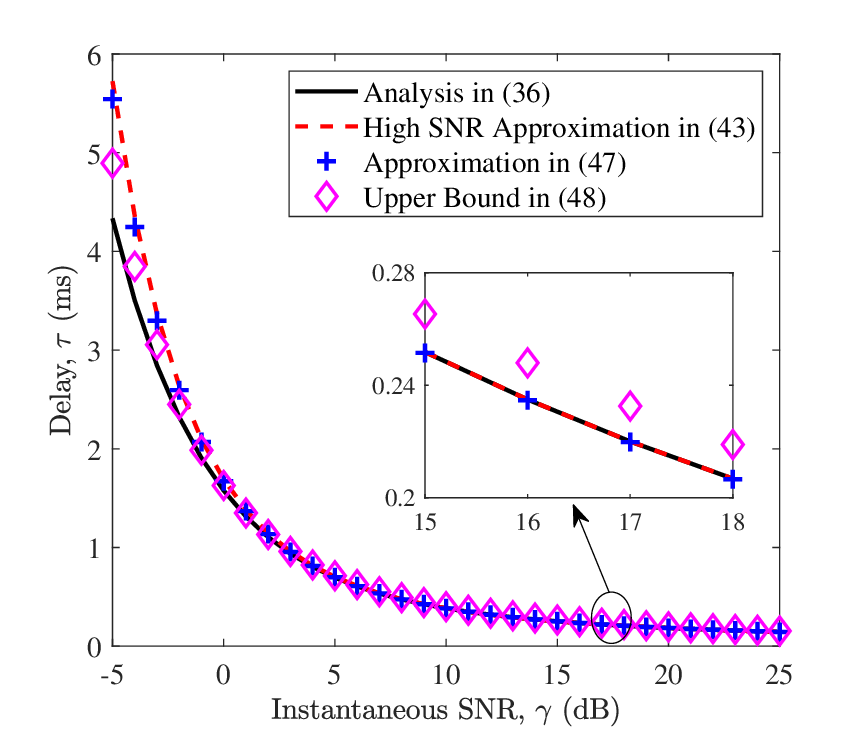} 
	\caption{The tightness of the approximations in \eqref{eq43}, \eqref{eq47} and \eqref{eq48}.}      
	\label{fig.6}                      
\end{figure}

Figure \ref{fig.6} validates the accuracy of the approximation and upper bound of transmission delay in the FBL regime (i.e., Lemma \ref{lemma7}). It can be observed that at low SNR, the approximation in \eqref{eq43} serves as an upper bound on the theoretical delay. While at high SNR, the approximation in \eqref{eq43} shows an excellent fit with the theoretical curve, demonstrating the accuracy of \eqref{eq43}. Due to the accuracy of the binomial approximation, the approximation in \eqref{eq47} is consistent with the approximation in \eqref{eq43}. The upper bound in \eqref{eq48} is close to the theoretical delay curve throughout the SNR range, which indicates that the upper bound is a tight upper bound for the theoretical delay. Therefore, it is accurate to approximate the delay as \eqref{eq48}.

\begin{figure}[htbp] 
	\centering  
	\subfigure[]{
		\label{fig7.sub.1}
		\includegraphics[width=0.43\textwidth]{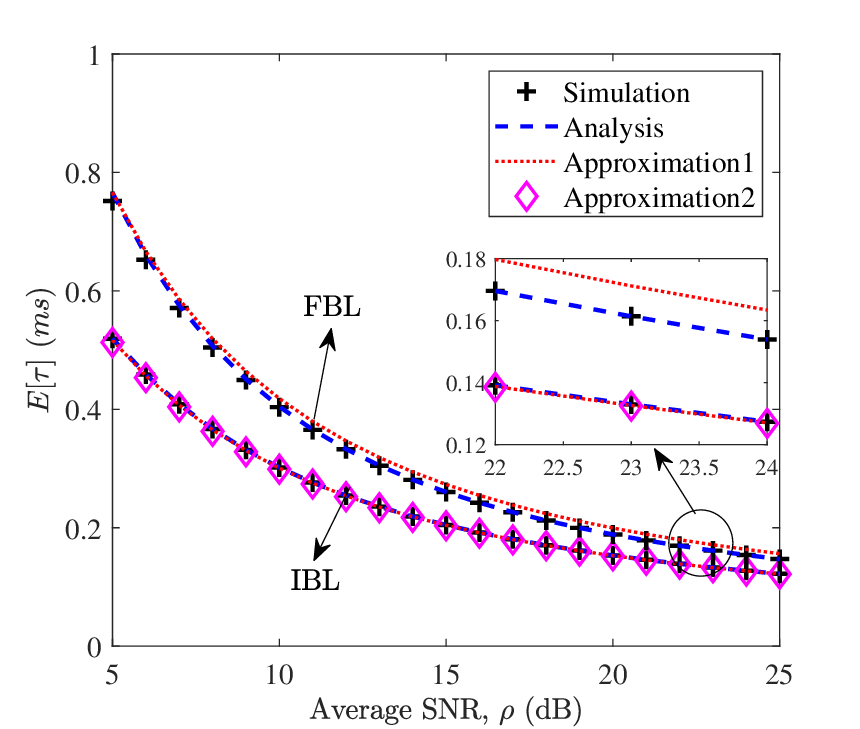}}	
	\subfigure[]{
		\label{fig7.sub.2}
		\includegraphics[width=0.43\textwidth]{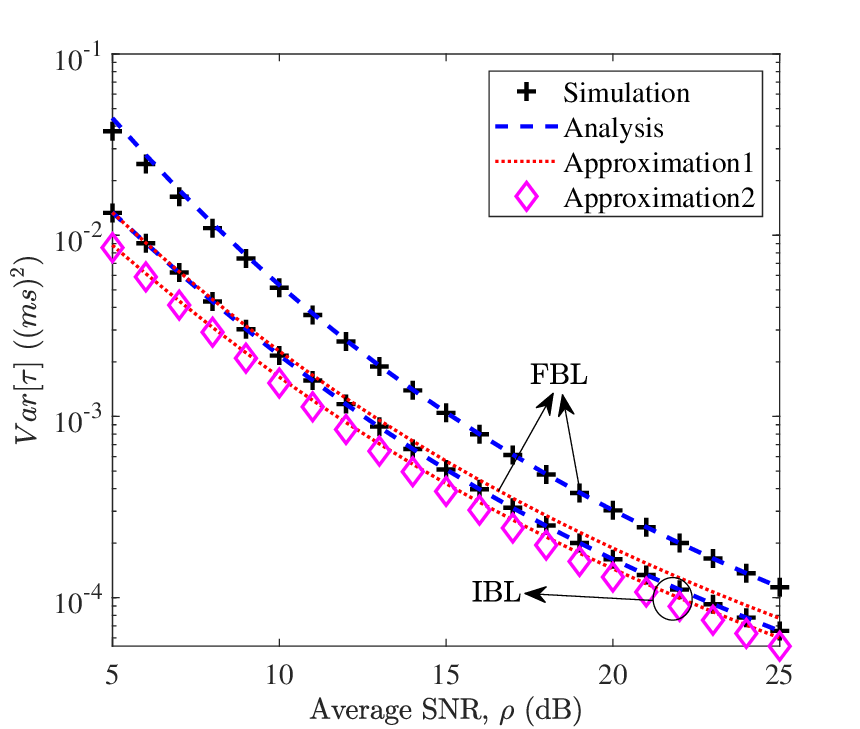}}
	\caption{The impact of average SNR on the delay performance. (a) Average delay. (b) Jitter.}
	\label{fig.7}
\end{figure}

\begin{figure}[htbp]  
	\centering  
	\subfigure[]{
		\label{fig8.sub.1}
		\includegraphics[width=0.43\textwidth]{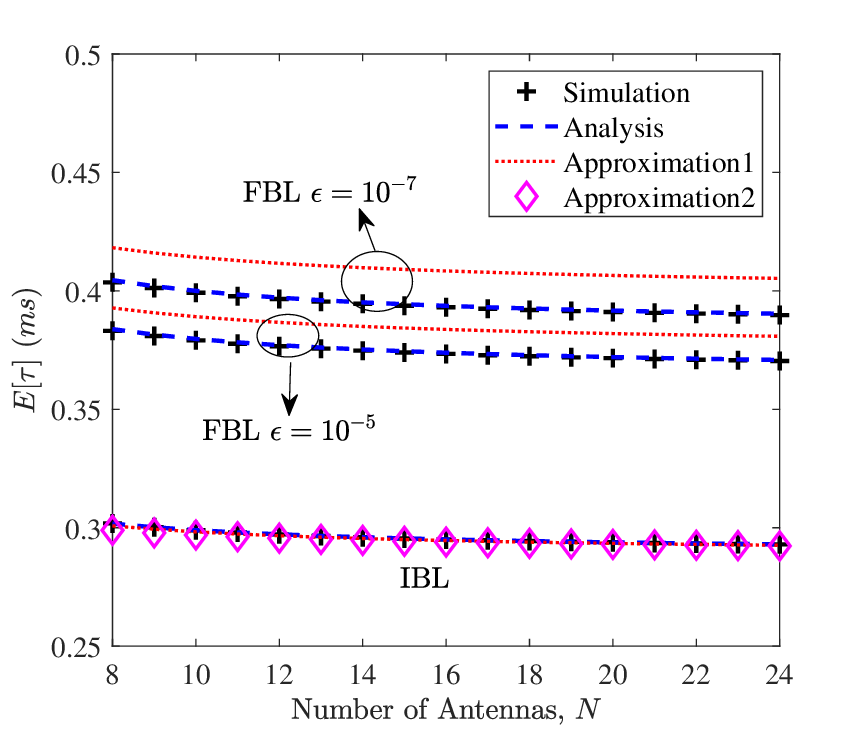}}	
	\subfigure[]{
		\label{fig8.sub.2}
		\includegraphics[width=0.43\textwidth]{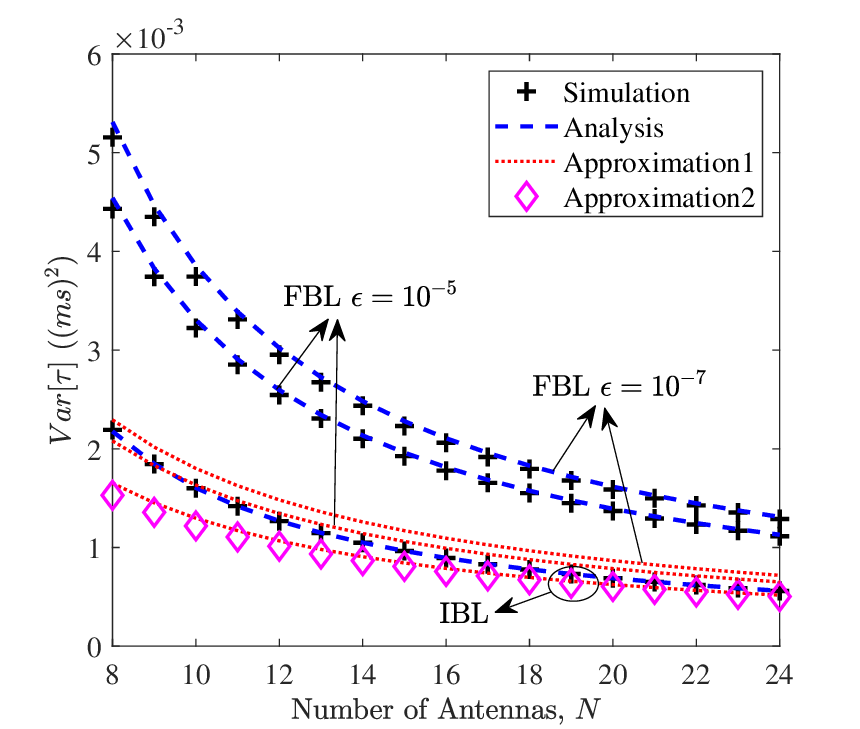}}
	\caption{The impact of the number of antennas on the delay performance. (a) Average delay. (b) Jitter.}
	\label{fig.8}
\end{figure}

To gain useful insights, we focus on the system parameters that affect the average delay and jitter. The effects of the average SNR and the number of antennas on the average delay and jitter are investigated in Figs. \ref{fig.7} and \ref{fig.8}, respectively. In Figs. \ref{fig.7} and \ref{fig.8}, the label ``Approximation1'' refers to the results obtained from Theorem \ref{theorem1} for IBL transmission and Theorem \ref{theorem3} for FBL transmission, respectively, while the label ``Approximation2'' represents the results obtained from Theorem \ref{theorem2} for IBL transmission. Figures \ref{fig.7} and \ref{fig.8} verify the accuracy of the expressions of average delay and jitter (i.e., Theorems \ref{theorem1}, \ref{theorem2}, and \ref{theorem3}). As shown in Fig. \ref{fig.7}, the mean and variance of the transmission delay decrease significantly with the average SNR, indicating that increasing the transmit power improves the delay performance. In the IBL regime, the four curves of the average delay almost overlap, and the variance curves are close to each other, which validates the accuracy of Theorem \ref{theorem1} and Theorem \ref{theorem2}. In the FBL regime, the analysis results are almost identical with the simulation results, but the approximate results deviate from the simulation results. It is mainly caused by the following two reasons: 1) an upper bound of delay is adopted in Theorem \ref{theorem3}; and 2) the mean and variance formulas \eqref{eq25} and \eqref{eq26} are derived by approximating the first two terms of the Taylor series, which may introduce some errors. In Fig. \ref{fig.8}, we can see that as the number of antennas increases, the average delay remains almost unchanged, while the variance decreases first and then remains unchanged. This indicates that the multi-antenna transmission is able to reduce jitter, but the performance no longer changes when the number of antennas is sufficiently large. Additionally, reducing the reliability requirement improves the delay performance. This suggests that latency and reliability are two conflicting performance metrics. Hence, trade-offs are inevitable in system design and performance optimization.

\subsection{Delay Violation Probability} \label{sec5.4}

\begin{figure}[htbp]     
	\centering       
	\includegraphics[width=0.43\textwidth]{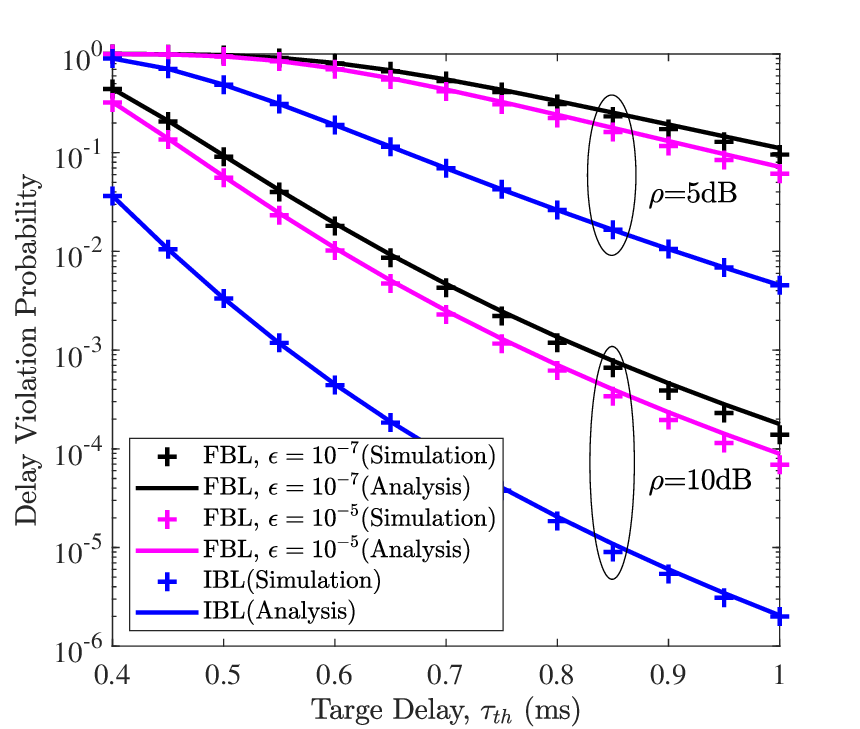} 
	\caption{The impact of target daley on the delay violation probability with different $\rho$ and $\epsilon$.}      
	\label{fig.9}                      
\end{figure}

\begin{figure}[htbp]      
	\centering       
	\includegraphics[width=0.43\textwidth]{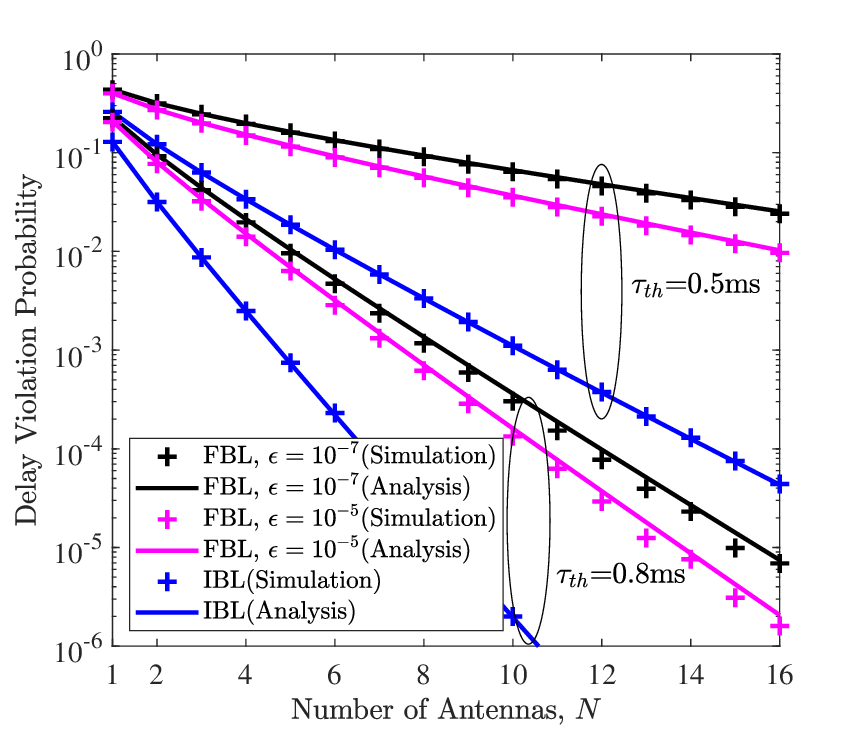} 
	\caption{The impact of the number of antennas on the delay violation probability with different $\tau_{th}$ and $\epsilon$.}      
	\label{fig.10}                      
\end{figure}

Figures \ref{fig.9} and \ref{fig.10} simulate the impact of system parameters on the delay violation probability. The number of channel realizations is $10^7$. As seen in Fig. \ref{fig.9}, the delay violation probability decreases with the target delay. Increasing the average SNR significantly reduces the delay violation probability for a given target delay. In the FBL regime, decreasing the reliability requirement decreases the delay violation probability. In addition, the delay violation probability calculated using Shannon capacity is significantly lower than that in the FBL regime, which indicates that the delay inscribed based on Shannon capacity is not accurate enough when the blocklength is short. Figure \ref{fig.10} depicts delay violation probability versus the number of antennas. It can be seen that the delay violation probability decreases with the number of antennas. 

\subsection{Design Insights} \label{sec5.5}
Based on the theoretical analysis and simulation results mentioned above, we have gained valuable insights into the design of delay-deterministic wireless networks.

\begin{itemize}
	\item  In the single-antenna scenario, the PDF of transmission delay has a large tail distribution. This is primarily due to the deep fading of the wireless channel and the limited transmission resources. Multi-antenna transmission achieves diversity gain by leveraging the independence between multi-antenna channels, which greatly reduces the probability of deep fading and shapes the distribution of transmission delay. Moreover, increasing the number of antennas reduces jitter. However, when the number of antennas becomes large, the impact on jitter becomes less noticeable. Therefore, multi-antenna transmit diversity can be used to improve delay performance. The number of antennas needs to be selected carefully, taking into account both hardware cost and delay performance. At the same time, beamforming can be considered to further reduce jitter.
	
	\item The delay performance in the FBL regime, including the statistical distribution, delay violation probability, average delay, and jitter, is significantly worse than that in the IBL regime. This is because the effects of noise and interference on the signal cannot be averaged out when the blocklength is short \cite{Durisi}. Therefore, for short packet communication scenarios such as industrial control, it is necessary to carry out research on delay performance analysis and optimization based on the finite blocklength transmission, which is still in its early stages.
	
	\item  In the FBL regime, decoding errors are unavoidable, and reliability and delay are conflicting performance metrics. Therefore, achieving deterministic low latency and high reliability requires careful trade-offs and rational design.
	
\end{itemize}
\section{Conclusion} \label{sec6}
This paper has analyzed the performance of transmission delay. In particular, we have considered both IBL and FBL transmissions. In the IBL regime, we have analyzed the distributions of transmission rate and transmission delay and derived the closed-form expressions for average delay, jitter, and delay violation probability. In the FBL regime, due to the complexity of the theoretical PDF expression of the transmission delay, we have approximated the delay at high SNR. Based on this approximation, the approximate delay distribution and delay violation probability have been obtained, and the closed-form expressions for average delay and jitter have been derived. The simulation results have validated the accuracy of the theoretical analysis and approximate calculations presented in this paper. The comparison of the results in the single-antenna and multi-antenna scenarios has illustrated that multi-antenna transmission can improve both the delay performance and the tail distribution of delay. The comparison of the results under the IBL and FBL regimes has demonstrated that the analytical results based on the traditional Shannon capacity overestimate the delay performance and are inaccurate. In addition, the system parameters that affect the average delay and jitter have been analyzed. Based on the theoretical analysis results presented in this paper, further investigation can be conducted on the resource allocation algorithm to achieve deterministic low latency communication. This aspect will be explored in our future research.

\appendices 
\section{Proof of the Theorem \ref{theorem1}}  \label{appendix1}
Assuming that the mean and variance of random variable $X$ are $\mu_X$ and $\sigma^2_X$, according to \cite{Benaroya}, the mean and variance of random variable function $f(X)$ can be approximated by the following equations:
\begin{equation*}  
	\mathbb{E} [f(X)] \approx f(\mu_X) + \frac{f''(\mu_X) }{2} \sigma_X^2, \tag{53} \label{eq53}
\end{equation*}
\begin{equation*} 
	\text{Var} [f(X)] \approx [f'(\mu_X)]^2 \sigma_X^2 - \frac{[f''(\mu_X)]^2}{4} \sigma_X^4. \tag{54} \label{eq54}
\end{equation*}
Considering $R^{-1}$ as a function of the random variable $R$, according to \eqref{eq53} and \eqref{eq54}, we have
\begin{equation*}   
	\mathbb{E} [R^{-1}] \approx \frac{1}{ \mathbb{E} [R] } + \frac{\text{Var}[R]}{(\mathbb{E}[R])^3}	= \frac{\mathbb{E}[R^2]}{\mathbb{E}[R]^3} , \tag{55} \label{eq55}
\end{equation*}
\begin{align*} 
	\text{Var}  [R^{-1}] &\approx \frac{\text{Var}[R]}{ (\mathbb{E}[R])^4 } - \frac{ (\text{Var}[R])^2 }{(\mathbb{E}[R])^6}\\ &= -\frac{\mathbb{E}[R^2]^2}{\mathbb{E}[R]^6} + \frac{3\mathbb{E}[R^2]}{\mathbb{E}[R]^4} - \frac{2}{\mathbb{E}[R]^2}. \tag{56} \label{eq56}
\end{align*}
From \eqref{eq55} and \eqref{eq56}, it can be seen that the closed-form expressions for the average delay and jitter require the first-order and second-order moments of the transmission rate. Note that the approximate formulas \eqref{eq53} and \eqref{eq54} are also used in \cite{R. Hashemi} and \cite{Jung}.

The moment of the transmission rate is
\begin{equation*}  
	\mathbb{E} [R^i ] = \int_0^{\infty} \left[\log_2(1+x)\right]^i f_{\gamma}(x) dx. \tag{57} \label{eq57}
\end{equation*}
In \eqref{eq57}, the main challenge arises from the logarithmic function. We can approximate the logarithmic function with
\begin{equation*} 
	\ln t \approx b t^{\frac{1}{b}} - b, \tag{58} \label{eq58} 
\end{equation*}
where $b$ is a large constant, and the approximation is sufficiently accurate when $b$ is large. The accuracy of the approximation in \eqref{eq58} will be proved by simulations in Section \ref{sec5}. 

Based on \eqref{eq58}, the first-order and second-order moments of the transmission rate are given by \eqref{eq59} and \eqref{eq60}, respectively, as shown at the bottom of the next page. In \eqref{eq59}, (a) is achieved by setting $ t = x+1 $, (b) is obtained by the binomial formula \cite[eq. (1.111)]{Gradshteyn}, (c) is based on the approximation in \eqref{eq58}, and (d) is based on \cite[eq. (3.381.3)]{Gradshteyn}. Equation \eqref{eq60} is calculated similarly to \eqref{eq59}.

By substituting \eqref{eq55}, \eqref{eq56}, \eqref{eq59}, and \eqref{eq60} into \eqref{eq23} and \eqref{eq24}, we obtain the closed-form expressions for the average delay and jitter. 

\begin{figure*}[hb]
	\hrulefill   
	\begin{align*} \mathbb{E} [R] & = \frac{1}{\ln2} \int_{0}^{\infty}  \ln \left( 1+x\right) \frac{N^N x^{N-1}  e^{-\frac{N}{\rho}x} }{\rho^{N}\Gamma(N)} dx \\
		&\overset{(a)}{=} \frac{N^N e^{\frac{N}{\rho}}}{\ln2 \rho^{N}\Gamma(N)} \int_{1}^{\infty} (t-1)^{N-1}  e^{-\frac{N}{\rho}t} \ln t dt \\
		& \overset{(b)}{=} \frac{N^N e^{\frac{N}{\rho}}}{\ln2 \rho^{N}\Gamma(N)} \sum_{k=0}^{N-1} \begin{pmatrix} N-1\\k \end{pmatrix} (-1)^{N-1-k} \int_{1}^{\infty} t^k e^{-\frac{N}{\rho}t} \ln t dt \\
		& \overset{(c)}{\approx } \frac{N^N e^{\frac{N}{\rho}}}{\ln2 \rho^{N}\Gamma(N)} \sum_{k=0}^{N-1} \begin{pmatrix} N-1\\k \end{pmatrix} (-1)^{N-1-k} 
		\left[ b\int_{1}^{\infty} t^{\frac{1}{b}+k} e^{-\frac{N}{\rho}t} dx - b\int_{1}^{\infty} t^k e^{-\frac{N}{\rho}t} dt \right] \\
		& \overset{(d)}{=} \frac{N^N e^{\frac{N}{\rho}}}{\ln2 \rho^{N}\Gamma(N)} \sum_{k=0}^{N-1} \begin{pmatrix} N-1\\k \end{pmatrix} (-1)^{N-1-k} 
		\left[ b\left(\frac{\rho}{N}\right)^{ \frac{1}{b}+k+1 }\Gamma\left(\frac{1}{b}+k+1,\frac{N}{\rho}\right) - b\left(\frac{\rho}{N}\right)^{k+1}\Gamma\left(k+1,\frac{N}{\rho}\right) \right] 
		\tag{59} \label{eq59}
	\end{align*} 
\end{figure*}

\begin{figure*}[hb]	
	\begin{align*}  
		\mathbb{E} [R^2]  =& \frac{1}{(\ln2)^2}\int_{0}^{\infty} \left[\ln \left( 1+x\right)\right]^2   \frac{N^N x^{N-1}  e^{-\frac{N}{\rho}x} }{\rho^{N}\Gamma(N)} dx  \\
		\approx & \frac{N^N e^{\frac{N}{\rho}}}{(\ln2)^2\rho^{N}\Gamma(N)} \sum_{k=0}^{N-1} \begin{pmatrix} N-1\\k \end{pmatrix} (-1)^{N-1-k}  \left[ b^2\int_{1}^{\infty} t^{\frac{2}{b}+k} e^{-\frac{N}{\rho}t} dt - 2b^2\int_{1}^{\infty} t^{\frac{1}{b}+k} e^{-\frac{N}{\rho}t} dt + b^2\int_{1}^{\infty} t^k e^{-\frac{N}{\rho}t} dt \right] \\
		= & \frac{N^N e^{\frac{N}{\rho}}}{(\ln2)^2\rho^{N}\Gamma(N)} \sum_{k=0}^{N-1} \begin{pmatrix} N-1\\k \end{pmatrix} (-1)^{N-1-k} \\ & \times \left[ b^2\left(\frac{\rho}{N}\right)^{ \frac{2}{b}+k+1 }\Gamma\left(\frac{2}{b}+k+1,\frac{N}{\rho}\right) - 2b^2\left(\frac{\rho}{N}\right)^{\frac{1}{b}+k+1 }\Gamma\left(\frac{1}{b}+k+1,\frac{N}{\rho}\right)  +b^2\left(\frac{\rho}{N}\right)^{k+1}\Gamma\left(k+1,\frac{N}{\rho}\right) \right] \tag{60} \label{eq60}
	\end{align*} 
\end{figure*}

\section{Proof of the Lemma \ref{lemma5}} \label{appendix2}
The CDF of the transmission delay is given by
\begin{equation*}  
	F_{\tau}(t) = \Pr \left\{ n(\gamma,L,\epsilon) \le Bt \right\}. \tag{61} \label{eq61}
\end{equation*}
According to the implicit function derivation theorem, for any given $L$ and $\epsilon$, the blocklength $n$ is a strictly monotonic decreasing function of $\gamma$ \cite{S. Xu}, and \eqref{eq61} can be rewritten as
\begin{equation*}  
	F_{\tau}(t) = \Pr \left\{ \gamma \ge n^{-1} \left( Bt \right)  \right\}  
	= 1- F_{\gamma}(n^{-1} \left( Bt \right)), \tag{62} \label{eq62}
\end{equation*}
where $\gamma=n^{-1}(Bt)$ is the solution of $  F(Bt,\gamma) = 0 $, which satisfies
\begin{equation*} 
	\ln \left( 1+\gamma\right) - \sqrt{1-\left( 1+\gamma\right)^{-2}} \frac{Q^{-1}\left( \epsilon \right)}{\sqrt{Bt}} - \frac{L\ln2}{Bt} = 0. \tag{63} \label{eq63}
\end{equation*}
Since it is difficult to solve \eqref{eq63} directly, considering the binomial approximation  $(1+z)^{\alpha} \approx 1 + \alpha z$, $|z|<1$, $|\alpha z|<1$ \cite{R. Hashemi}, we have
\begin{equation*} 
	\sqrt{1- \frac{1}{(1+\gamma)^2}}  \approx 1 - \frac{1}{2(1+\gamma)^2}. \tag{64} \label{eq64}
\end{equation*}
The accuracy of the approximation in \eqref{eq64} will be proved by simulations in Section \ref{sec5}. Based on \eqref{eq64}, \eqref{eq63} can be transformed as
\begin{equation*} 
	\ln \left[ e^{ -\frac{Q^{-1}\left( \epsilon \right)}{\sqrt{Bt}} -\frac{L\ln2}{Bt}} \left( 1+\gamma\right) \right] = - \frac{Q^{-1}\left( \epsilon \right)}{2\sqrt{Bt}} \frac{1}{\left( 1+\gamma\right)^{2}}. \tag{65} \label{eq65}
\end{equation*}
Let $\theta = \frac{Q^{-1}\left( \epsilon \right)}{\sqrt{Bt}} $, $
\phi = e^{ -\frac{Q^{-1}\left( \epsilon \right)}{\sqrt{Bt}} - \frac{L\ln2}{Bt}} $, and $
\eta = \ln \left[\phi\left( 1+\gamma\right) \right] $, and then we have $ 1+\gamma = \phi^{-1} e^{\eta} $. \eqref{eq65} can be rewritten as
\begin{equation*} 
	\eta = -\frac{\theta}{2} \phi^2 e^{-2\eta}. \tag{66} \label{eq66} 
\end{equation*}
Let $\kappa = 2\eta$ and rewrite \eqref{eq66} as
\begin{equation*}  \kappa e^{\kappa} = -\theta \phi^2, \tag{67} \label{eq67} 
\end{equation*}
where $\kappa$ can be expressed by the Lambert W function, i.e., $\kappa = \mathcal{W} \left( -\theta \phi^2 \right) $ \cite{Olver}. Thus, we have
\begin{equation*}  
	n^{-1} \left( Bt \right) = \phi^{-1} e^{\frac{1}{2}\mathcal{W} \left( -\theta \phi^2 \right) } - 1. \tag{68} \label{eq68} 
\end{equation*}
According to the Taylor series expansion of the Lambert W function, $ \mathcal{W} \left( x \right) = \sum_{m=1}^{\infty} \frac{(-m)^{m-1}}{m!} x^m $ \cite{Olver}, we have
\begin{equation*}  
	n^{-1} \left( Bt \right) = \phi^{-1} e^{-\frac{1}{2} \sum_{m=1}^{\infty} \frac{m^{m-1}}{m!}\theta^m \phi^{2m} } - 1. \tag{69} \label{eq69} 
\end{equation*}
By substituting \eqref{eq69} and \eqref{eq5} into \eqref{eq62}, we obtain the CDF of the transmission delay, and derivation of the CDF yields the PDF. 

Note that the references \cite{S. He} and \cite{Huang} have also solved \eqref{eq63}. However, the references \cite{S. He} and \cite{Huang} solved \eqref{eq63} directly and expressed the results in terms of the general Lambert W function, while we applied the approximation in \eqref{eq64} and expressed the results in terms of the Lambert W function. In contrast, the computational results presented in this paper are simpler.

\section{Proof of the Theorem \ref{theorem3}} \label{appendix3}
Let $ \Phi = \frac{1}{\ln(1+\gamma)} $, according to \eqref{eq53} and \eqref{eq54}, and based on the upper bound of the delay in \eqref{eq48}, the mean and variance of the delay are approximated as follows
\begin{align*}  
	\mathbb{E}\left[ \tau \right] \approx& 
	\frac{ 2L\ln 2+\left(Q^{-1}\left( \epsilon \right) \right)^2 }{ 2B } \mathbb{E}\left[ \Phi \right] \\&
	+ \frac{\sqrt{L\ln2} Q^{-1}\left( \epsilon \right)}{B} \mathbb{E}\left[ \Phi^{\frac{3}{2}} \right] , \tag{70} \label{eq70}
\end{align*}
\begin{align*}  
	\text{Var} \left[ \tau \right] \approx &
	\left[\frac{ 2L\ln 2+\left(Q^{-1}\left( \epsilon \right) \right)^2 }{ 2B }\right]^2 \text{Var} \left[ \Phi \right] \\&
	+ \frac{L\ln2 \left(Q^{-1}\left( \epsilon \right)\right)^2}{B^2} \text{Var} \left[ \Phi^{\frac{3}{2}} \right], \tag{71} \label{eq71}
\end{align*}
where
\begin{equation*}  
	\mathbb{E}\left[\Phi\right] = \frac{1}{\ln2} \left[ \frac{1}{ \mathbb{E} [R] } + \frac{\text{Var}[R]}{(\mathbb{E}[R])^3} \right], \tag{72} \label{eq72}
\end{equation*}
\begin{equation*} 
	\text{Var}\left[\Phi\right] = \frac{1}{(\ln2)^2} \left[ \frac{\text{Var}[R]}{ (\mathbb{E}[R])^4 } - \frac{ (\text{Var}[R])^2 }{(\mathbb{E}[R])^6} \right], \tag{73} \label{eq73}
\end{equation*}
\begin{equation*}  
	\mathbb{E}[\Phi^{\frac{3}{2}}] = \frac{1}{(\ln2)^{\frac{3}{2}}} \left[ \frac{1}{ (\mathbb{E}[R])^{\frac{3}{2}} } + \frac{15\text{Var}[R]}{8(\mathbb{E}[R])^{\frac{7}{2}} } \right], \tag{74} \label{eq74} 
\end{equation*}
\begin{equation*} 
	\text{Var} [\Phi^{\frac{3}{2}}] = \frac{1}{\left(\ln2\right)^3} \left[\frac{9 \text{Var} [R] }{4\left(\mathbb{E}[R]\right)^5} 
	- \frac{225 \left(\text{Var} [R]\right)^2 }{64\left(\mathbb{E}[R]\right)^7} \right]. \tag{75} \label{eq75}
\end{equation*}
Substituting \eqref{eq72}–\eqref{eq75} into \eqref{eq70} and \eqref{eq71} yields the closed-form expression for the mean delay and jitter in \eqref{eq51} and \eqref{eq52}.

\end{document}